%% file: Dsept2603AG.tex
\documentclass[12pt,letterpaper]{amsart}
\usepackage{amsmath,amsthm,amsfonts,latexsym,epsfig,eepic, epic}
\input xy
\xyoption{all}
\newtheorem{theorem}{Theorem}[section]
\newtheorem{conjecture}[theorem]{Conjecture}
\newtheorem{lemma}[theorem]{Lemma}
\newtheorem{proposition}[theorem]{Proposition}

   \newlength{\baseunit}               
   \newcount{\numlines}                
\newcommand{\bpf}{\noindent {\em Proof.  }}
\newcommand{\epf}{\qed \vspace{+10pt}}

\newcommand{\Sg}{\mathfrak{S}}

\newtheorem{corollary}[theorem]{Corollary}

\newcommand{\Z}{\mathbb{Z}}
\newcommand{\bA}{\mathbf{A}}
\newcommand{\bB}{\mathbf{B}}
\newcommand{\bC}{\mathbf{C}}
\newcommand{\bD}{\mathbf{D}}
\newcommand{\HH}{\mathbf{H}}
\newcommand{\hh}{\mathbf{h}}
\newcommand{\B}{S}

\newcommand{\R}{\mathbb{R}}
\newcommand{\C}{\mathbb{C}}

\newcommand{\proj}{\mathbb P}

\newcommand{\cP}{{\mathcal{P}}}
\newcommand{\cm}{{\mathcal{M}}}
\newcommand{\cmbar}{\overline{\cm}}

\newcommand{\al}{\alpha}

\newcommand{\be}{\beta}

\newcommand{\de}{\delta}
\newcommand{\De}{\Delta}
\newcommand{\si}{\sigma}
\newcommand{\la}{\lambda}
\newcommand{\Aut}{\operatorname{Aut}}
\newcommand{\Autal}{\left| \Aut \al \right|}
\newcommand{\Autbe}{\left| \Aut \be \right|}
\newcommand{\Autalbe}{\left| \Aut \al \right| \left| \Aut \be \right|}
\newcommand{\Pic}{\operatorname{Pic}}
\newcommand{\oPic}{\overline{\Pic}}

\newcommand{\genus}[1]{#1}
\newcommand{\tH}{\tilde{H}}
\newcommand{\llangle}{\pmb{\langle} \pmb{\langle}}
\newcommand{\rrangle}{\pmb{\rangle} \pmb{\rangle}}
\newcommand{\subst}[2]{  {\substack { {#1} \\ {#2}}} } 

\newcommand{\cited}{}

\newcommand{\lremind}[1]{{}}
\newcommand{\cut}[1]{}

\newif\ifstanford 
\newif\ifwaterloo 
\stanfordtrue 
\waterloofalse 
\newcommand{\stanford}[1]{\ifstanford {#1} \else {} \fi}
\newcommand{\waterloo}[1]{\ifwaterloo {#1} \else {} \fi}

\begin{document}
\pagestyle{plain}
\title{{\large {Towards the geometry of double Hurwitz numbers}}
}
\author{I. P. Goulden,
D. M. Jackson
and
R. Vakil}
\address{Department of Combinatorics and
Optimization, University of Waterloo}
\email{ipgoulden@math.uwaterloo.ca}
\address{Department of Combinatorics and
Optimization, University of Waterloo}
\email{dmjackson@math.uwaterloo.ca}
\address{Department of Mathematics, Stanford University}
\email{vakil@math.stanford.edu}
\thanks{The first two authors are partially supported by
NSERC grants. 
The third author is partially supported by NSF grant
DMS--0228011, NSF CAREER grant DMS--0238532, and an Alfred P. Sloan
Research Fellowship.
}
\date{Friday, September 26, 2003}

\begin{abstract}
  Double Hurwitz numbers count branched covers of $\mathbb{C} \proj^1$
  with fixed branch points, with simple branching required over all
  but two points $0$ and $\infty$, and the branching over $0$ and
  $\infty$ points specified by partitions of the degree (with $m$ and
  $n$ parts respectively).  Single Hurwitz numbers (or more usually,
  Hurwitz numbers) have a rich structure, explored by many authors in
  fields as diverse as algebraic geometry, symplectic geometry,
  combinatorics, representation theory, and mathematical physics.  A
  remarkable formula of Ekedahl, Lando, M.  Shapiro, and Vainshtein
  (the ELSV formula) relates single Hurwitz numbers to intersection
  theory on the moduli space of curves.  This connection has led to
  many consequences, including Okounkov and Pandharipande's proof
  of Witten's conjecture (Kontsevich's theorem) connecting
  intersection theory on the moduli space of curves to integrable
  systems.
  
  In this paper, we determine the structure of double Hurwitz numbers
  using techniques from geometry, algebra, and representation theory.
  Our motivation is geometric: we give strong evidence that double
  Hurwitz numbers are top intersections on a moduli space of curves
  with a line bundle (a universal Picard variety).  In particular, we
  prove a piecewise-polynomiality result analogous to that implied by
  the ELSV formula.
  In the case $m=1$ (complete branching over one point) and
  $n$ is arbitrary, we conjecture an ELSV-type formula, and show it to be
  true in genus $0$ and $1$.  The corresponding Witten-type
  correlation function has a better structure than that for single
  Hurwitz numbers, and we show that it satisfies many geometric
  properties, such as the string and dilaton equations, and a genus
  expansion ansatz analogous to that of Itzykson and Zuber.
  We give a symmetric function description of the
  double Hurwitz generating series, which leads to explicit formulae
  for double Hurwitz numbers with given $m$ and $n$, as a function of
  genus.   In the
  case where $m$ is fixed but not necessarily $1$, 
  we prove a topological recursion on the corresponding
  generating series, which leads to closed-form expressions for double
  Hurwitz numbers and an analogue of the Goulden-Jackson polynomiality
  conjecture (an early conjectural variant of the ELSV formula).
  
  Coupled with earlier work connecting double Hurwitz numbers to
  integrable systems, an ELSV-type formula would translate all of this
  structure on double Hurwitz numbers to the intersection theory of
  the universal Picard variety described earlier.  
\end{abstract}
\maketitle
\tableofcontents

\section{Introduction}
If $\alpha=(\al_1,\cdots,\al_m)$ and $ \beta=(\be_1,\cdots,\be_n)$ are 
partitions of a positive integer $d$, the {\em double Hurwitz number}
$H^\genus{g}_{\al, \be}$ is the number of genus $g$ branched covers 
of $\mathbb{C} \proj^1$ with branching corresponding to $\al$ and $\be$
over $0$ and $\infty$ respectively, and an appropriate number $r$ 
of other
fixed simple branched points (determined by the Riemann-Hurwitz formula,
and depending on $g$, $m$, $n$).
For simple branching, the monodromy of the sheets is a transposition. 
To simplify the exposition, we assume that the points
mapping to $0$ and $\infty$ are {\em labelled}.  Thus the double Hurwitz
numbers under this convention are $\left| \Aut \al \right| \left| \Aut \be \right|$ larger
than they would be under the convention in~\cite{gjvk}.

Our (long-term) goal
is to understand the structure of double Hurwitz numbers, and in
particular to determine the possible form of an ELSV-type formula
expressing double Hurwitz numbers in terms of intersection theory on
some universal Picard variety, that we surmise is close to the one
defined by Caporaso in \cite{lucia}.  (M. Shapiro has made significant progress in
determining what this space might be \cite{spc}.)

\subsection{Motivation from single Hurwitz numbers}
\label{elsvsub}\lremind{elsvsub}Our 
methods are extensions of the combinatorial and
character-theoretic methods that we have used in the very
well-developed theory of {\em single Hurwitz numbers}
$H^\mathbf{g}_\alpha$, where all but possibly one branch point have
simple branching.  (They are usually called ``Hurwitz numbers,'' but
we have added the term ``single'' to distinguish them from the double
Hurwitz numbers.)  Single Hurwitz numbers have proved to have
surprising connections to geometry, including the moduli space of
curves.  (For a remarkable recent link to the Hilbert scheme of points
on a surface, see for example \cite[p.~2]{surface} and
\cite{vasserot}.)  Our intent is to draw similar connections in the
case of the double Hurwitz numbers.  We wish to use the
representation-theoretic and combinatorial structure of double
Hurwitz numbers to understand the intersection theory of a conjectural 
universal
Picard variety, in analogy with the connection between single Hurwitz
numbers and the moduli space of curves, as shown in the following
diagram.

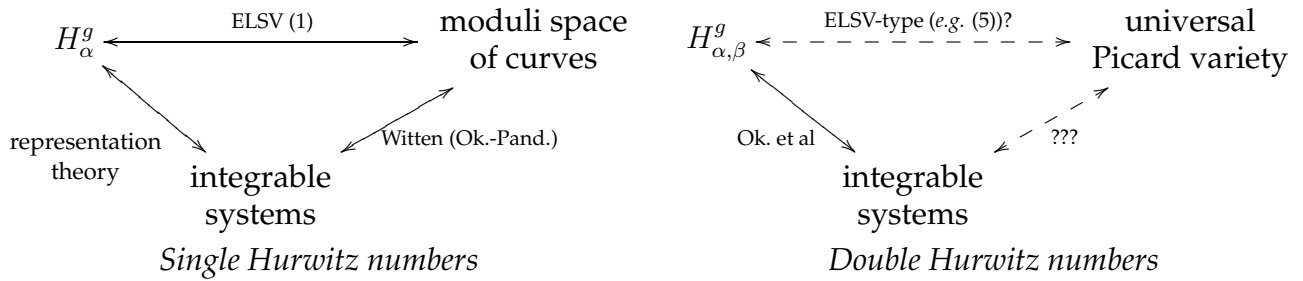
\stanford{$$
\begin{array}{cc}
\xymatrix{H^{\genus{g}}_{\al} \ar@{<->}[0,2]^{\text{ELSV \eqref{elsv}} \quad \quad}  
 \ar@{<->}[dr]_{\scriptsize {\begin{array}{c} \text{representation} \\ \text{theory}\end{array}}} & &  {\begin{array}{c}\text{moduli space} \\ \text{of curves} \end{array}}
\ar@{<->}[dl]^{\text{\quad Witten (Ok.-Pand.)}} 
\\
& 
{\begin{array}{c} \text{integrable} \\ \text{systems}\end{array}}
}
&
\xymatrix{H^{\genus{g}}_{\al, \be} \ar@{<-->}[0,2]^{\text{ELSV-type ({\em e.g.} \eqref{earlyconj})? \quad \quad \quad }}  
 \ar@{<->}[dr]_{\text{Ok.\ et al}} & & {\begin{array}{c}\text{universal} \\ \text{Picard variety} \end{array}}
\ar@{<-->}[dl]^{\text{???}} 
\\
&
{ \begin{array}{c} \text{integrable} \\ \text{systems} \end{array}}
}
\\
\text{\em Single Hurwitz numbers} & \text{\em Double Hurwitz numbers} \\
\end{array}
$$
Understanding this would give, for example, Toda constraints on the topology of the
moduli of curves.}
\waterloo{Understanding this would give, for example, Toda constraints on the topology of the
moduli of curves.}

The history of single Hurwitz numbers is too
long to elaborate here (and our bibliography omits many foundational
articles), but we wish to draw the reader's attention 
to ideas leading, in particular, to the ELSV formula.

The ELSV formula (\cite{elsv1, elsv2}, see also \cite{gv}) asserts that:
\begin{equation} \label{elsv}
\boxed{H^\genus{g}_{\al} = C(g,\alpha)
    \int_{\cmbar_{g,m}} \frac { 1 -
    \la_1 + \la_2 \cdots \pm \la_g} { ( 1- \al_1 \psi_1) \cdots (1 -
    \al_m \psi_m)}}
\end{equation}
where
\begin{eqnarray}\label{eCg}
    C(g,\alpha) = r! \prod_{i=1}^m \frac {\al_i^{\al_i}} {\al_i!} 
\end{eqnarray}
is a  scaling factor (it clearly has a combinatorial interpretation).
Here $\cmbar_{g,m}$ is  Deligne and Mumford's compactification of the
moduli space of genus $g$  curves with $m$ marked points,
$\la_k$ is a certain codimension $k$ class, and
$\psi_i$ is a certain codimension $1$ class.  We refer the reader to the original papers
for precise definitions, which we will not need.  (The original ELSV
formula includes a factor of $\left| \Aut \al \right|$ in the
denominator, but as stated earlier, we are considering the points over
$\infty$, or equivalently the parts of $\al$, to be labelled.)  The
right hand side should be interpreted by expanding the integrand formally,
and capping the terms of degree $\dim \cmbar_{g,m}=3g-3+m$ with
the fundamental class 
$[\cmbar_{g,m}]$.

The ELSV formula \eqref{elsv} implies that
\begin{eqnarray}\label{eHCP}
H^{\genus{g}}_{\al} = C(g,\al) P^{\genus{g}}_m (\al_1, \dots, \al_m),
\end{eqnarray}
where $P^{\genus{g}}_m$ is a
{\em polynomial whose terms have total degrees between $2g-3+m$ and $3g-3+m=\dim \cmbar_{g,m}$}.  
The coefficients of this polynomial are precisely the
top intersections on the moduli space of curves involving $\psi$-classes
and up to one $\la$-class, often  written, using Witten's notation, 
as:
\lremind{wittennotation}
\begin{equation}
\label{wittennotation}
\boxed{\langle \tau_{a_1} \dots \tau_{a_m} \la_k \rangle_g
:= \int_{\cmbar_{g,m}} \psi_1^{a_1} \cdots \psi_m^{a_m} \la_k 
= (-1)^k \left[ \al_1^{a_1} \cdots \al_m^{a_m}  \right] P^\genus{g}_m(\al_1, \dots, \al_m)}
\end{equation}
when $\sum a_i + k = 3g-3+m$, and $0$ otherwise. (Here we use the
notation $[A]B$ for the coefficient of $A$ in $B$.)
This ELSV polynomiality is related to (and implies, by \cite[Thm.\ 3.2]{gjvk}) 
an earlier conjecture of Goulden and Jackson, describing
the form of the generating series for single Hurwitz numbers of
genus $g$ (\cite[Conj.\ 1.2]{gj2}, see also \cite[Conj.\ 1.4]{gjvn}).
The conjecture asserts that after a change of variables, the single Hurwitz
generating series is ``polynomial'' (in the sense that its scaled coefficients
are polynomials). The Goulden-Jackson conjecture
is in fact a genus expansion ansatz for Hurwitz numbers analogous
to the ansatz of Itzykson and Zuber (\cite[(5.32)]{iz}, proved in 
\cite{eyy, gjvk}).  
ELSV polynomiality is related to Goulden-Jackson polynomiality, an
early conjectural variant of the ELSV formula, by
a change of variables arising from Lagrange inversion \cite[Thm.\ 2.5]{gjvk}.

Hence, in developing the theory of double Hurwitz numbers, we seek
some sort of polynomiality (in this case, piecewise polynomiality)
that will tell us something about the moduli space in the background
(such as its dimension), as well as a Goulden-Jackson-type conjecture
or ansatz.

\subsection{Summary of results}
In Section~\ref{homogeneity}, we use ribbon graphs to establish that
double Hurwitz numbers (with fixed $m$ and $n$) are piecewise
polynomial of degree up to $4g-3+m+n$ (Piecewise Polynomiality Thm.\ 
\ref{homthm}), with {\em no} scaling factor analogous to $C(g,\al)$.
More precisely, for fixed $m$ and $n$, we show that
$H^\genus{g}_{(\al_1, \dots, \al_m), (\be_1, \dots, \be_n)}$ 
counts the number of lattice points in certain polytopes, and as
the $\al_i$ and $\be_j$ vary, the facets move.
Further, we conjecture that the degree is bounded below by $2g-3+m+n$
(Conj.\ ~\ref{stronghom}), and verify this conjecture in genus $0$,
and for $m$ or $n=1$ (Cor.\ \ref{stronghomtrue}).  We give an example
($(g,m,n)=(0,2,2)$) showing that it is not polynomial in general.

In Section~\ref{m1}, we consider the 
case $m=1$ (``one-part double Hurwitz numbers''), which corresponds to
double Hurwitz numbers with complete branching over $0$.  One-part double
Hurwitz numbers have a particularly tractable structure.
In
particular, they {\em are} polynomial: for fixed $g$, $n$,
$H^\genus{g}_{(d), (\be_1, \dots, \be_n)}$ is a polynomial
in $\be_1$, \dots, $\be_n$.
Theorem~\ref{onepartg} gives two formulae for these
numbers (one in terms of the series $\sinh x / x$ and the other 
an explicit expression) generalizing
formulae of both Shapiro-Shapiro-Vainshtein \cite[Thm.\ 6]{ssv} and
Goulden-Jackson \cite[Thm.\ ~3.2]{gjcacti}.
As an application, we prove
polynomiality, and in particular show that the resulting polynomials
have simple expressions in terms of character theory. 
Based on this polynomiality, we conjecture an
ELSV-type formula for one-part double Hurwitz numbers (Conj.\ 
~\ref{onepartconj}):\lremind{earlyconj}
\begin{equation}
\label{earlyconj}
\boxed{H^{\genus{g}}_{(d), \be} = r!  d 
\int_{\oPic_{g,n}}  \frac  { \Lambda_0 - \Lambda_2 + \cdots \pm \Lambda_{2g}}
{( 1 - \beta_1 \psi_1) \cdots (1 - \be_n \psi_n)}
}
\end{equation}
Here, $\Pic_{g,n}$ is the moduli space
of smooth genus $g$ curves with $n$ distinct labelled
smooth points together with a line bundle; for example, the points of
$\Pic_{g,n}$ correspond to ordered triples (smooth genus $g$ curve $C$, $n$ distinct labeled
points on $C$, line bundle on $C$).
The space $\oPic_{g,n}$ is some as-yet-undetermined compactification of 
$\Pic_{g,n}$ , supporting classes $\psi_i$ and
$\Lambda_{2k}$, satisfying properties described in Conj.\ \ref{onepartconj}.
As with the ELSV formula \eqref{elsv}, the right side of \eqref{earlyconj} should be interpreted by
expanding the integrand formally, and capping the terms
of dimension $4g-3+n$ with $[\oPic_{g,n}]$. 
The most speculative part of this conjecture is the identification
of the $(4g-3+n)$-dimensional moduli space with a compactification 
of $\Pic_{g,n}$ (see the Remarks following Conj.\ \ref{onepartconj}).

Motivated by this conjecture, we define a symbol $\llangle \cdot
\rrangle_g$, the 
 analogue of $\langle \cdot \rangle_g$,
 by the first equality of
\begin{eqnarray}\label{eWPL}
\boxed{\llangle \tau_{b_1} \cdots \tau_{b_n} \Lambda_{2k} \rrangle_g
:=
(-1)^{k} \left[ \be_1^{b_1} \cdots \be_n^{b_n} \right] \left(
\frac{  H^{\genus{g}}_{(d), \be }}
{r! d }\right) 
= \int_{\oPic_{g,n}} \psi_1^{b_1} \cdots \psi_n^{b_n} \Lambda_{2k},}
\end{eqnarray}
so that Conjecture~\ref{onepartconj} (or \eqref{earlyconj})
would imply the second equality.  (All parts of \eqref{eWPL} are
zero unless $\sum b_i +2k = 4g-3+n$.
Also, we note that the definition of $\llangle \cdot \rrangle_g$ is
independent of the conjecture.) We show that this symbol satisfies
many properties analogous to those proved by Faber and Pandharipande
for $\langle \cdot \rangle_g$, including integrals over
$\cmbar_{g,1}$, and the $\la_g$-theorem; we generalize these further.
We then prove a genus expansion ansatz for $\llangle \cdot \rrangle_g$
in the style of Itzykson-Zuber (\cite{iz}~Thm.\ \ref{ResultI}).  As
consequences, we prove that $\llangle \cdot \rrangle_g$ satisfies the
string and dilaton equations, and verify the ELSV-type conjecture in
genus $0$ and $1$.  A proof of Conjecture~\ref{onepartconj} would
translate all of this structure associated with double Hurwitz numbers
to the intersection theory of the universal Picard variety.

In Section~\ref{character2}, we give a simple formula for the double
Hurwitz generating series in terms of Schur symmetric functions.  As an
application, we give explicit formulae for double Hurwitz numbers
$H^\genus{g}_{\al, \be}$ for fixed $\al$ and $\be$, in terms of
linear combinations of $g$th powers of prescribed integers,
extending work of Kuleshov and M. Shapiro \cite{ks}.
Although this section is placed after Section~\ref{m1}, it can be
read independently of Section~\ref{m1}.

In Section~\ref{marb}, we consider $m${\em-part Hurwitz numbers} (those
with $m=l(\al)$ fixed and $\be$ arbitrary).  As remarked earlier,
polynomiality fails in this case in general, but we still find strong
suggestions of geometric structure.  We define a (symmetrized)
generating series $\HH^\genus{g}_m$ for these numbers, and show that it
satisfies a topological recursion (in $g$, $m$)
(Theorems~\ref{gensymm}, ~\ref{eqminfour}, ~\ref{intgenusg}).  The
existence of such a recursion is somewhat surprising as, unlike other
known recursions in Gromov-Witten theory (involving the geometry of
the {source} curve), it is not a low-genus phenomenon. (The one
exception is the Toda recursion of \cite{ptoda, o}, which also
deals with double Hurwitz numbers.)  We
use this recursion to derive closed expressions for $\HH^{\genus{g}}_m$
for small $(g,m)$, and to conjecture a general form (Conj.\ 
\ref{conjecture1}), in analogy with the original Goulden-Jackson
polynomiality conjecture of \cite{gj2}.

\noindent {\bf Acknowledgments.}
We are grateful to R. Hain, A. Knutson, 
A. Okounkov and M. Shapiro for helpful discussions,
and to M. Shapiro for sharing his ideas.
The third author thanks 
R. Pandharipande for explaining the ribbon
graph construction for single Hurwitz numbers, which is (very mildly)
generalized in Section~\ref{homogeneity}.
We are grateful to I. Dolgachev and A. Barvinok for pointing out 
the  reference \cite{mc}.  J. Bryan and E. Miller suggested
improvements to the manuscript.

%

\subsection{Earlier evidence of structure in double Hurwitz numbers}
Our work is motivated by several recent suggestions of strong
structure of double Hurwitz numbers.  Most strikingly, Okounkov proved
that the generating series $H$ for double Hurwitz numbers is a
$\tau$-function for the Toda hierarchy of Ueno and Takasaki \cite{o},
in the course of resolving a conjecture of Pandharipande's on {single}
Hurwitz numbers \cite{ptoda}; see also their joint work \cite{op1,op2,op3}.
Dijkgraaf's earlier description \cite{dijkgraaf} of 
Hurwitz numbers
where the target has genus $1$ and all branching is simple, and his
unexpected discovery that the corresponding generating series is
essentially a quasi-modular form, is also suggestive, as such Hurwitz
numbers can be written (by means of a generalized join-cut equation)
in terms of double Hurwitz numbers (where $\al=\be$).  This
quasi-modularity was generalized by Bloch and Okounkov \cite{bo}.

Signs of structure for fixed $g$ (and fixed number of points) provides
a clue to the existence of a connection between double Hurwitz
numbers and the moduli of curves (with additional structure), and even
suggests the form of the connection, as was the case for single
Hurwitz numbers.  Evidence for this comes from recent work of Lando
and D. Zvonkine \cite{lz}, Kuleshov and M. Shapiro \cite{ks}, and others.

We note that double Hurwitz numbers are relative Gromov-Witten
invariants (see for example \cite{li1} in the algebraic category, and
earlier definitions in the symplectic category \cite{lr, ip}), and
hence are necessarily top intersections on a moduli space (of
relative stable maps).  Techniques of Okounkov and Pandharipande
\cite{op1, op2, op3} can be used to study double Hurwitz numbers in
this guise. A second promising approach, relating more general Hurwitz
numbers to intersections on moduli spaces of curves, is due to Shadrin
\cite{shadrin} building on work of Ionel \cite{ionel}.  We expect that
some of our results are probably obtainable by one of these two
approaches.  However, we were unable to use them to prove any of the
conjectures and, in particular, we could prove no ELSV-type formula.

We also alert the reader to important recent work on
Hurwitz numbers due to Kuleshov, Lando, M. Shapiro, and D. Zvonkine
\cite{l, lz, ks, z}. 

\subsection{Notation and background} \label{notation} 
\lremind{notation}Throughout, 
the partitions $\al$ and $\be$ have $m$ and $n$ parts, respectively. We
use $l(\alpha)$ for the number of parts of $\alpha$, and $\left| \alpha \right|$
for the sum of the parts of $\alpha$. If $\left| \alpha \right|=d$, we say
$\alpha$ is a partition of $d$, and write $\alpha\vdash d$. 
For a partition $\alpha =(\alpha_1,\ldots )$, let $\Aut \al$ be the group of
permutations of $\{ 1, \dots, l(\al) \}$ fixing $(\al_1, \dots, \al_{l(\al)})$. 
Hence, if $\alpha$ has $a_i$ parts equal to $i$, $i\geq 1$, then 
$\left| \Aut \alpha  \right| =\prod_{i\geq 1} a_i!$.  For indeterminates $p_1,\ldots$ and
$q_1,\ldots$, we write $p_{\alpha}= \prod_{i \geq 1}p_{\alpha_i}$ and
$q_{\alpha}=\prod_{i\geq 1}q_{\alpha_i}$. Let ${\mathcal{C}}_{\alpha}$ denote the
conjugacy class of the symmetric group $\Sg_d$ indexed by $\alpha$,
so $\left| \mathcal{C}_\al
\right| =d! /\Autal \prod_i \al_i$.
We use the notation $[A]B$ for
 the coefficient of monomial $A$ in a formal power series series $B$.

Genus will in general be denoted by superscript.
Let \lremind{rh2}
\begin{equation}r^{\genus{g}}_{\al,\be} := -2+2g+m+n. \label{rh2}
\end{equation} 
When the context permits, we shall abbreviate this to $r.$

A summary of other globally defined notation is in the table below.

\begin{center}
\begin{tabular}{|l|l|} \hline
$\langle \cdot \rangle_g$, $\tau_i$
    & Witten symbol \eqref{wittennotation} \\[5pt]  
$Q$, $w$, $w_i$, $\mu$, $\mu_i$, $Q_i$
    & Lagrange's Implicit Function Theorem ~\ref{lagthm} \\[5pt]   
\rule{0pt}{5pt}
$H^\genus{g}_{\al, \be}$, $\tH^\genus{g}_{\al, \be}$, $H$, $\tH$
    & double Hurwitz numbers and series, Section~\ref{dHnnot} \\[5pt] 
$\Theta_m$, $\HH^\genus{g}_m$, $\HH^\genus{g}_{m,i}$
    & symmetrization operator, symmetrized genus $g$ $m$-part \\[5pt]
    & \quad \quad Hurwitz function, and its derivatives Sect.\ \ref{symop} \\[5pt] 
$P^\genus{g}_{m,n}$ 
    & Piecewise Polynomiality theorem ~\ref{homthm} \\[5pt]  
$\mathsf{E}^{\alpha}$, $\mathsf{K}_{\alpha}$
    & character theory \eqref{char1}, \eqref{char3} \\[5pt]  
$s_{\theta}, h_i, p_i$ 
    & symmetric functions (Schur, complete, power sum), \\[5pt]
    &  Sections~\ref{m45},  ~\ref{character2} \\[5pt]  
$N_i$, $c_i = N_i -\de_{i,1}$, $S_{2j}$
    & functions of $\be$, Section~\ref{m11}  \\[5pt] 
$B_{2k}$, $\xi_{2k}$ (and $\xi_{2 \la}$), $f_{2k}$, $v_{2k}$
    & coefficients of $\frac x {e^x-1} + x/2$ (Bernoulli),
      $\log \frac {\sinh x} x$ (Thm.\ \ref{onepartg}), \\[5pt]
    & \quad \quad $\frac {x/2} { \sinh (x/2)},\;\;\frac{2}{x}\sinh (x/2)$ 
      (Thm.~\ref{Wittansx}), 
       \\[5pt]  
$\llangle \cdot \rrangle_g$, $\oPic_{g,n}$, $\Lambda_{2k}$
    & ELSV-type Conjecture ~\ref{onepartconj} \\[5pt]  
$Q^{(i)}(t)$, $Q^{(\la )}(t)$
    & genus expansion ansatz Theorem~\ref{ResultI} \\[5pt]  
$\hh_m^\genus{g}= \Gamma \HH^\genus{g}_m$, $\hh^\genus{g}_{m,i}$  
    & transform of $\HH^\genus{g}_m$, and its partial derivatives, Sect.\ \ref{ss:Grss} \\[5pt]  \hline
\end{tabular}
\end{center}

\subsubsection{Double Hurwitz numbers.} \label{dHnnot} 
\lremind{dHnnot}As 
described earlier, let the {\em double Hurwitz number} $H^{\genus{g}}_{\al,\be}$
denote the number of degree $d$ branched covers of $\C \proj^1$ by a
genus $g$ (connected) Riemann surface, with $r+2$ branch
points, of which  $r=r^{\genus{g}}_{\al,\be}$ are
simple,  and two ($0$ and $\infty$, say) have 
branching given by $\al$ and $\be$, respectively.  Then
\eqref{rh2} is the Riemann-Hurwitz formula.
If a cover has automorphism group $G$, it is counted with multiplicity
$1/\left|G\right| $.  For example, $H^{\genus{0}}_{(d),(d)} = 1/d$.  The points above
$0$ and $\infty$ are taken to be {\em unlabelled}.

The {\em possibly disconnected} double Hurwitz numbers $\tH^{\genus{g}}_{\al, \be}$ are
defined in the same way except the covers are not required to be connected.

The double Hurwitz numbers may be characterized in terms of the symmetric group through
the monodromy of the sheets around the branch points.
This axiomatization is essentially due to Hurwitz~\cite{hurwitz};
the proof relies on the Riemann existence theorem.
\begin{proposition}[Hurwitz axioms] 
\label{huraxioms}\lremind{huraxioms}For $\alpha ,\beta \vdash d$, 
$H^{\genus{g}}_{\alpha ,\beta}$
is equal to $\Autalbe /d!$ times the number
of $(\sigma ,\tau_1,\ldots ,\tau_r,\gamma )$, such
that
\begin{enumerate}
\item[H1.] $\sigma\in{\mathcal{C}}_{\beta},\;\; \gamma\in
{\mathcal{C}}_{\alpha},
\;\; \tau_1,\ldots ,\tau_r$ are transpositions on $\{ 1,\ldots ,d\}$,
\item[H2.] $\tau_r\cdots \tau_1\sigma = \gamma$,
\item[H3.] $r=r^{\genus{g}}_{\al,\be}$, and
\item[H4.] the group generated by $\sigma ,\tau_1,\ldots ,\tau_r$ acts transitively
on $\{ 1,\ldots ,d\}$.
\end{enumerate}
The number $\tH^{\genus{g}}_{\al, \be}$ is equal to $\Autalbe /d!$ times
the number of $(\sigma ,\tau_1,\ldots ,\tau_r,\gamma )$ satisfying H1--H3.
\end{proposition}
If  $(\sigma ,\tau_1,\ldots ,\tau_r)$
satisfies H1--H3, we call it an {\em ordered factorization}
of $\gamma$, and if it also satisfies H4, we call it
a {\em transitive ordered factorization}.

The {\em double Hurwitz (generating) series} $H$ 
for double Hurwitz numbers is given by
\lremind{defH}
\begin{equation} \label{defH}
H = \sum_{g\geq 0,d\geq1}
\sum_{\alpha ,\beta\vdash d}
y^g z^d
p_{\alpha}q_{\beta} u^{l(\beta )}
\frac{H^{\genus{g}}_{\alpha,\beta}}{r^{\genus{g}}_{\al,\be}! 
\Autalbe},
\end{equation}
and $\tH$ is the analogous generating series for the possibly
disconnected double Hurwitz numbers.  Then 
$\tH = e^H$,
by a general enumerative result (see, {\em e.g.,}~\cite[Lem.\ 3.2.16]{gjbook}).
(The earliest reference we know for this result is, appropriately
enough, in work of Hurwitz.)

The following result is obtained by using the axiomatization above,
and by studying the effect that multiplication by a final
transposition has on the cutting and joining of cycles in the cycle
decomposition of the product of the remaining factors.  The details of
the proof are essentially the same as that of ~\cite[Lem.\ 2.2]{gj0} and
~\cite[Lem.\ 3.1]{gjvn}, and are therefore suppressed.  A
geometric proof involves pinching a loop separating the
target $\C \proj^1$ into two disks, one of which contains only one simple
branch point and the branch point corresponding to $\beta$.

\begin{lemma}[Join-cut equation] \label{jc}
\begin{equation}\label{joincut}
\left( \sum_{i\geq 1}p_i \frac{\partial}{\partial p_i}
+u\frac{\partial}{\partial u}+2y\frac{\partial}{\partial y}
-2\right) H=\tfrac{1}{2}\sum_{i,j\geq 1} \left(
ijp_{i+j}\frac{\partial H}{\partial p_i}
\frac{\partial H}{\partial p_j}
+(i+j)p_i p_j \frac{\partial H}{\partial p_{i+j}}
+ijp_{i+j} y \frac{{\partial}^2H}{\partial p_i
\partial p_j}  \right)
\end{equation}
with initial conditions $\left[ z^i p_i q_i u \right]H =\frac{1}{i}$ for $i\geq 1$.
\end{lemma}

Substituting $u\frac{\partial}{\partial u}H= \sum_{i\geq 1}q_i
\frac{\partial}{\partial q_i}H$ yields the usual, more symmetric
version.  But the above formulation will be more convenient for our
purposes.


\subsubsection{The symmetrization operator $\Theta_m$,
and the symmetrized double Hurwitz generating series $\HH^{\genus{g}}_m$}
\label{symop} \lremind{symop}The 
linear symmetrization operator $\Theta_m$ is defined by
\begin{eqnarray}\label{elso}
\Theta _m(p_{\alpha})=\sum_{\sigma\in \Sg_m}x_{\sigma (1)}^{\alpha _1}
\cdots x_{\sigma (m)}^{\alpha _m}
\end{eqnarray}
if $l(\alpha )=m$, and  zero otherwise. (It is {\em not} a ring homomorphism.)
The properties of $\Theta_m$ we require appear as Lemmas~4.1, 4.2, 4.3 in~\cite{gjvn}.
Note that $\Theta _m(p_{\alpha})$ has a close relationship with the monomial
symmetric function $m_\alpha$ since
$$\Theta _m(p_{\alpha}) = \Autal  m_{\alpha}(x_1,\ldots ,x_m).$$
We shall study in detail the
symmetrization $\sum_{m\geq 1, g\geq 0} \HH^{\genus{g}}_m y^g$, of $H$
where
\begin{eqnarray} \nonumber
\HH^{\genus{g}}_m(x_1,\ldots ,x_m) &=& \left[ y^g \right]\left.\Theta_m \left(
H\right)\right|_{z=1}, \\
\label{symtransf}
&=&\sum_{d\geq 1}
\sum_{\alpha,\beta\vdash d \atop l(\al)=m}
m_{\alpha}(x_1,\ldots ,x_m)q_{\beta}u^{l(\beta)} 
   \frac{H^{\genus{g}}_{\alpha ,\beta}}{r^{\genus{g}}_{\al,\be}! \Autbe},
\end{eqnarray}
 \lremind{symtransf}for $m\ge1, g\geq 0$.

%
In other words, the redundant variable $z$ is eliminated, and
$\HH^{\genus{g}}_m$ is a generating series containing information about
genus $g$ double Hurwitz numbers (where $\al$ has $m$ parts).

We use the notation
$$\HH^{\genus{g}}_{j,i}=x_i \frac {\partial \HH^{\genus{g}}_j} {\partial x_i}.$$

\subsubsection{Lagrange's Implicit Function Theorem.}
\label{lintro} \lremind{lintro}We shall 
make repeated use of the following form of 
Lagrange's Implicit Function Theorem,
(see, {\em e.g.,}~\cite[Sect.\ 1.2]{gjbook} for a proof) concerning
the solution of certain formal functional equations.

\begin{theorem}[Lagrange] \label{lagthm} \lremind{lagthm}Let 
$\phi(\la)$ be an invertible formal power series in an
indeterminate $\lambda.$ Then the
functional equation
$$v=x\phi (v)$$
has a unique formal power series solution $v=v(x).$
Moreover, if $f$ is a formal power series, then
\begin{eqnarray}
\label{lagthm1}
f(v(x)) &=& f(0)+\sum_{n\geq
1}\frac{x^n}{n}\left[ \lambda^{n-1} \right]\frac{d f(\lambda )}{d\lambda}
\phi (\lambda )^n \quad \quad \mbox{and}
\\
\label{lagthm2}
\frac{f(v(x))}{v}\frac{xdv(x)}{dx} &=& \sum_{n\geq 0}x^n \left[ \lambda
^n \right]f(\lambda )
\phi (\lambda )^n. \end{eqnarray}
\end{theorem}

We apply Lagrange's theorem to the functional
equation \lremind{funceqw}
\begin{equation}\label{funceqw} w=x e^{uQ(w)}, \end{equation}
where $$Q(t)=\sum_{j\geq 1} q_j t^j,$$
the series in the indeterminates $q_j$ that record the parts of $\be$ in the double
Hurwitz series \eqref{defH}.
The following observations and notation will be used extensively.
By differentiating the functional
equation with respect to $x$, and $u$, we obtain
\begin{equation}\label{lagpartial}
x\frac{\partial w}{\partial x}=w\mu(w),\;\;\;\;
\frac{\partial w}{\partial u}=wQ(w)\mu(w),\;\;\;\;
\mbox{where $\mu(t)=\frac{1}{1-utQ^{\prime}(t)}$,}
\end{equation}
and we therefore have the operator identity
\begin{equation}\label{opident}
\frac{x\partial}{\partial x}=\mu(w)\frac{w\partial}{\partial w}.
\end{equation}
We shall use the notation $w_i=w(x_i)$, $\mu_i=\mu (w_i)$, and
$Q_i=Q(w_i)$, for $i=1,\ldots ,m$.


\section{Piecewise polynomiality} \label{homogeneity} 
\lremind{homogeneity}By 
analogy with the ELSV formula \eqref{elsv}, we consider double Hurwitz numbers
for fixed $g$, $m$, $n$ as functions in the parts of $\al$ and $\be$:
$$ P^{\genus{g}}_{m,n}(\al_1, \dots, \al_m, \be_1, \dots, \be_n) =
H^{\genus{g}}_{\al,\be} .$$
Here the domain is the set of
$(m+n)$-tuples of positive integers, where the sum of the first $m$
terms equals the sum of the remaining $n$. In contrast
with the single Hurwitz number case, the double Hurwitz numbers have
no scaling factor $C(g, \al, \be)$ 
(see ~(\ref{eHCP})).\lremind{homthm}

\begin{theorem}[Piecewise polynomiality]  
For fixed $m$, $n$, the double Hurwitz function
  $H^{\genus{g}}_{\al,\be} = P^{\genus{g}}_{m,n}$ is piecewise
  polynomial (in the parts of $\al$ and $\be$) of degrees up to
  $4g-3+m+n$.  The ``leading'' term of degree $4g-3+m+n$ is non-zero.
\label{homthm}\lremind{homthm}
\end{theorem}

By {\em non-zero leading term}, we mean that for fixed $\al$ and
$\be$, $P^{\genus{g}}_{m,n}(\al_1 t, \dots, \be_n t)$ (considered as a
function of $t \in \Z^+$) is a polynomial of degree $4g-3+m+n$.  In
fact, this leading term can be interpreted as the volume of a certain
polytope.
For example, 
$P^{\genus{0}}_{2,2}(\al_1, \al_2, \be_1, \be_2) = 2
\max(\al_1, \al_2, \be_1, \be_2)$, which has degree one.
(This can be shown by a straightforward
calculation,  either directly, or using Section~\ref{hompf},
or Cor.\ \ref{twotwo}.
See Cor.\ \ref{twotwo} for a calculation of $P^\genus{g}_{2,2}$
in general.)
In particular, unlike the case  of single Hurwitz numbers
(see~(\ref{eHCP})), $P^{\genus{g}}_{m,n}$ is not polynomial in 
general.

We conjecture further:
\lremind{stronghom}

\begin{conjecture}[Strong piecewise polynomiality] \label{stronghom}
  $P^{\genus{g}}_{m,n}$ is piecewise polynomial, with degrees between
  $2g-3+m+n$ and $4g-3+m+n$ inclusive.
\end{conjecture}

It is straightforward to check that
Conjecture~\ref{stronghom}
is true in genus $0$ (by induction using the joint-cut
equation \eqref{joincut}), and we will show that it holds when $m$ or
$n$ is $1$ (Cor.\ \ref{stronghomtrue}).  It may be possible to verify the
conjecture by refining the proof of the Piecewise Polynomiality
Theorem, but we were unable to do so.

\subsection{Proof of the Piecewise Polynomiality Theorem \ref{homthm}}\label{hompf}
We spend the rest of this section proving Theorem~\ref{homthm}.
Our strategy is to interpret double Hurwitz numbers as counting
lattice points in certain polytopes.
We use a combinatorial interpretation of double Hurwitz numbers that is a
straightforward extension of the
interpretation of single Hurwitz numbers given in \cite[Sect.\ 3.1.1]{op1}
(which is there shown to be equivalent to earlier graph interpretations
of \cite{arnold, ssv}).  The case $r=0$
is trivially verified (the double Hurwitz number is $1/d$ if $\al =
\be = (d)$ and $0$ otherwise), so we assume $r>0$.

Consider a branched cover of $\C \proj^1$ by a genus $g$ Riemann surface $S$,
with branching over $0$ and $\infty$ given by $\al$ and $\be$, and $r$
other branch points (as in Sect.\ \ref{notation}).  We may assume that
the $r$ branch points lie on the equator of the $\C \proj^1$, say at the
$r$ roots of unity.  Number the $r$ branch points $1$ through $r$ in
counterclockwise order around $0$.

We construct a ribbon graph on $S$ as follows.  The vertices are the $m$
preimages of $0$, denoted by $v_1$, \dots, $v_m$, where $v_i$
corresponds to $\al_i$.  For each of the $r$
branch points on the equator of the target $b_1$, \dots, $b_r$, consider the $d$
preimages of the geodesic (or radius) joining $b_r$ to $0$.  Two of
them meet the
corresponding ramification point on the source.  Together, they form
an edge joining two (possibly identical) of the vertices.  The
resulting graph on the genus $g$ surface has $m$ (labelled) vertices
and $r$ (labelled) edges.  There are $n$ (labelled) faces, each
homotopic to an open disk.  The faces correspond to the parts of $\be$:
each preimage of $\infty$ lies in a distinct face.  Call such a
structure a {\em labelled} (ribbon) {graph}.  (Euler's formula
$m-r+n=2-2g$ is equivalent to the Riemann-Hurwitz
formula \eqref{rh2}.)

Define a {\em corner} of this labelled graph to be the data consisting of a
vertex, two edges incident to the vertex and adjacent to each other
around the vertex, and the face between them (see Figure~\ref{cornereg}).

\stanford{
\begin{figure}[ht]
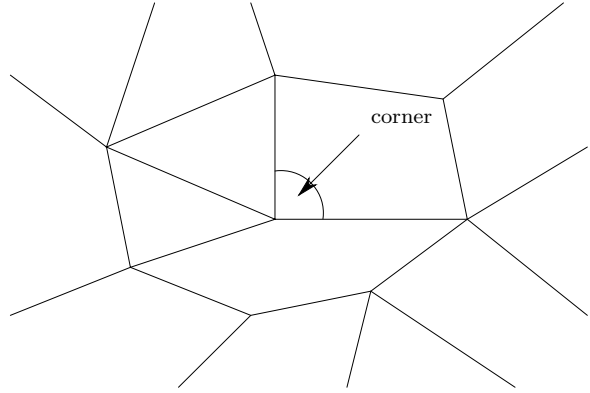

\begin{center}
\include{cornereg4}
\end{center}
\caption{An example of a corner in a fragment of a map  \label{cornereg} \lremind{cornereg}}
  \end{figure}
}
\waterloo{
 \begin{figure}[h]
 \begin{center}
  \includegraphics[scale=1.0]{cornereg4.eps}%
 \end{center}
 \caption{An example of a corner in a fragment of a map
      \label{cornereg} \lremind{cornereg}}
 \end{figure}
}

Now place a dot near $0$ on the target $\C \proj^1$, between the geodesics to
the branch points $r$ and $1$.  Place  dots on the source surface
$S$ at the $d$ preimages of the dot on the target $\C \proj^1$.

Then the number of dots near vertex $v_i$ is $\al_i$: a small circle
around $v_i$ maps to a loop winding $\al_i$ times around $0$.  Moreover,
any corner where edge $i$ is counterclockwise of $j$ and $i<j$ must
contain a dot.  Call such a corner a {\em descending corner}. The number of
dots in face $f_j$ is $\be_j$: move the  dot on the target
(together with its $d$ preimages) along a line of longitude until it is
near the pole $\infty$ (the $d$ preimages clearly do not cross any
edges en route), and repeat the earlier argument.

Thus each cover counted in the double Hurwitz number corresponds to a
combinatorial object: a labelled graph (with $m$ vertices, $r$ edges,
and $n$ faces, hence genus $g$), and a non-negative integer (number of
dots) associated to each corner, which is positive if the corner is
descending, such that the sum of the integers around vertex $i$ is
$\al_i$, and the sum of the integers in face $j$ is $\be_j$.

It is straightforward to check that the converse is true (using the
Riemann existence theorem, see for example \cite{arnold}): given such
a combinatorial structure, one gives the target sphere a complex
structure (with branch points at roots of unity), and this induces a
complex structure on the source surface.

Hence the number of double Hurwitz numbers is a sum over the set of
labelled graphs (with $m$ vertices, $r$ edges, and $n$ faces).  The
contribution of each labelled graph is the number of ways of assigning
non-negative numbers to each corner so that each descending corner is
assigned a positive integer, and such that the sum of numbers around
vertex $i$ is $\al_i$ and the sum of the integers in face $j$ is
$\be_j$.

For fixed $m$ and $n$, 
the contributions to $H^{\genus{g}}_{\al,\be}$ is the sum
over the same finite set of labelled graphs.  Hence to prove the
Piecewise Polynomiality Theorem it suffices to show that, for $\al$ and $\be$
fixed, the contribution of each labelled graph to $P^{\genus{g}}_{m,n}(t \al_1,
\dots, t \be_m)$ is a polynomial of degree $4g-3+m+n$.

This problem corresponds to counting points in a polytope as follows.
We have one variable for each corner (which is the corresponding
number of dots).  The number of corners is easily seen to be twice the
number of edges (see Fig.\ \ref{corneredge}), so there are $2r$
variables $z_1,\ldots,z_{2r}.$  We have one linear equation for each vertex (the sum of
the variables corresponding to corners incident to vertex $i$ must be
$\al_i$) and one for each face (the sum of the variables corresponding
to corners incident to face $j$ must be $\be_j$).  These equations are
dependent since the sum of the $m$ vertex relations is the sum of the $n$
face relations, {\em i.e.}  $\sum^m_{i=1} \al_i = \sum^n_{j=1} \be_j$.

\stanford{
\begin{figure}[ht]
\begin{center}
\include{corneredge2}
\end{center}
\caption{Twice the number of corners is four times the number of edges
\lremind{corneredge}\label{corneredge}}
\end{figure}
}
\waterloo{
 \begin{figure}[h]
 \begin{center}
  \includegraphics[scale=1.0]{corneredge.eps}%
 \end{center}
  \caption{Twice the number of corners is four times the number of edges
  \label{cornereg} \lremind{cornereg}}
 \end{figure}
}

There are no other dependencies,  {\em i.e.} the rank of the system is
$m+n-1$: suppose otherwise, that one of the equations, for example the
equation $eq_1$ corresponding to vertex $i$, were a linear combination
others modulo the sum relation.  Pick a face $j$ incident to that
vertex.  Let $z$ be the variable corresponding to the corner between
vertex $i$ and face $j$.  Discard the equation $eq_2$ corresponding to
that face (which is redundant because of the sum relation).  Then $z$
appears only in equation $eq_1$, and hence $eq_1$ cannot be a linear
combination of the other equations.

Thus the contribution to the double Hurwitz number
$P^{\genus{g}}_{m,n}(\al_1, \dots, \be_n)$ from this labelled graph
$\Gamma$ is the number of lattice points in a polytope
$\cP_\Gamma(\al_1, \dots, \be_n)$ of dimension $2r-(m+n-1) = 4g-3+m+n$
in $\R^{2r}$, lying in the linear subspaces defined by
\begin{equation}
\label{redblue}
\sum_{\substack {      {\text{corner $k$ incident}}   \\
{\text{ to vertex $i$}}  }   } z_k =t \al_i \quad
\quad \text{and} \quad \quad \sum_{
\substack{     {\text{corner $k$ incident}}   \\
  {\text{ to face   $j$}}  }   }   z_k = t \be_j,
\end{equation}
bounded by inequalities of the form $z_k
\geq 0$ or $z_k > 0$ (depending on whether corner $k$ is descending or
not).  Let $P_{\Gamma}(  \al_1, \dots, \be_m)$
be this contribution.

\begin{lemma} The vertices of the polytope $\cP_\Gamma(\al, \be)$
are lattice points, {\em i.e.,} the polytope is integral.
\end{lemma}

\bpf Let $p \in \R^{2r}$ be a point of the polytope.  We show that if
$p$ is not a lattice point, then $p$ lies in the interior of a line
segment contained in $\cP_{\Gamma}(\al, \be)$, and hence is not a
vertex.  Construct an auxiliary graph, where the vertices correspond
to corners $i$ of $\Gamma$ such that $z_i(v) \notin \Z$.  The
edges come in two colors.  Red edges join any two distinct vertices
incident to a common vertex, and blue edges join any two distinct
vertices incident to a common face.  By the first (resp.\ second) equality in
\eqref{redblue}, each vertex is incident to a red (resp.\ blue) edge.
Thus we may find a cycle of distinct vertices $v_1=v_{2w+1}$, \dots, $v_{2w}$
such that $v_{2i-1}$ and $v_{2i}$ (resp.\ $v_{2i}$ and $v_{2i+1}$)
are joined by a red (resp.\ blue) edge:
choose any $x_1$, and then subsequently choose $x_2$, $x_3$, 
etc.\ (such that $x_i$ and $x_{i+1}$ is joined by an appropriately
coloured edge) until the first repetition:  $x_j=x_k$ ($j<k$).
If $k-j$ is even, take $v_i = x_{j+i}$ ($1 \leq i \leq k-j$),
and if $k-j$ is odd, take $v_i = x_{j+i}$ ($1  \leq  i < k-j$).
(If $v_1$ and $v_2$ are joined by a blue edge rather than a red edge,
then cyclically permute the $v_i$ by one.)

Then for $\left| \epsilon  \right| < \min (z_{v_i})$, the point $p(\epsilon)$
given by 
$$z_j(\epsilon)= \left\{
\begin{array}{ll} 
z_j &  \text{for $j \notin \{ v_1, \dots, v_{2w} \}$} \\
z_j+\epsilon & \text{for $j=v_{\text{even}}$} \\
z_j-\epsilon & \text{for $j=v_{\text{odd}}$} \\
\end{array}
\right.
$$
satisfies \eqref{redblue} and $z_j(t) \geq 0$, and hence
also lies in the polytope. \epf
  
As $\cP_{\Gamma}(\al, \be)$ is an integral polytope, by Ehrhart's theorem
\cite{eh}, for $t$ a positive integer, $P_{\Gamma}( t \al_1,
\dots, t \be_m)$ is a polynomial in $t$ of degree precisely $4g-3+m+n$
(the {\em Ehrhart polynomial} of the polytope), with leading
coefficient equal to the volume of $\cP_{\Gamma}(\al, \be)$. 
Finally, we recall
the following well-known result \cite{mc}.

\begin{theorem}
Consider the polytopes in $\R^L$ (with coordinates $z_1$, \dots, $z_L$)
defined by equalities 
$$\sum_{i=1}^L \mu_{ij} z_i = \eta_j  \quad \quad (1 \leq j \leq \si)
\quad \quad \quad \text{and} \quad \quad \quad
\sum_{i=1}^L \nu_{ij} z_i \leq 
\zeta_j \quad \quad (1 \leq j
\leq \tau)$$
as $\eta_j$ and $\zeta_j$ vary (with $\mu_{ij}$ and $\nu_{ij}$ fixed).
When the polytope is integral (for given $\eta_j$ and $\zeta_j$),
define $U(\eta_1, \dots, \eta_{\si}, \zeta_1, \dots, \zeta_{\tau})$ to be
the number of lattice points contained therein.
Then the function $U$ is piecewise polynomial on its domain, of
degree equal to the dimension of the polytope.
\end{theorem}
Thus
as $\al$ and
$\be$ vary, the function $P^\genus{g}_{m,n}$ is piecewise polynomial, 
concluding the proof of the Piecewise Polynomiality Theorem~\ref{homthm}.


\section{One-part double Hurwitz numbers ($\al = (d)$)}
\label{m1}\lremind{m1}We 
use character theory to completely describe double Hurwitz numbers
where $\al$ has one part, which leads to a conjectural
formula in terms of intersection theory on a moduli space.
The particular results that are needed from character theory
are to be found in~\cite{macd}.

\subsection{One-part double Hurwitz numbers through characters}\label{m11}
\lremind{m11}In 
the group algebra $\C\Sg_d$, let
   ${\mathsf{K}}_{\alpha}:=\sum_{\sigma\in {\mathcal C}_{\alpha}}\sigma.$
Then $\{ {\mathsf{K}}_{\alpha},\; \alpha\vdash d\}$, is a
basis for the centre, and
if $\chi^{\alpha}_{\lambda}$ is the character of the irreducible
representation of $\Sg_d$ indexed by ${\mathcal{C}}_{\alpha}$, evaluated
at any element of ${\mathcal{C}}_{\lambda}$, then \lremind{char1}
\begin{equation}\label{char1}
\mathsf{E}^{\alpha}=\frac{\chi_{(1^d)}^{\alpha}}{d!}\sum_{\lambda\vdash d}
\chi ^{\alpha}_{\lambda}{\mathsf{K}}_{\lambda},\;\;\;\; \alpha\vdash d,
\end{equation}
gives a basis of orthogonal idempotents.  The inverse relations 
are \lremind{char3}
\begin{equation}\label{char3}
{\mathsf{K}}_{\alpha}=\left|{\mathcal{C}}_{\alpha} \right| \sum_{\lambda\vdash d}
\frac{\chi_{\alpha}^{\lambda}}{\chi_{(1^d)}^{\lambda}}\mathsf{E}^{\lambda},
\;\;\;\; \alpha\vdash d.
\end{equation}

The following result gives an expression for the double Hurwitz
number $H^{\genus{g}}_{(d),\beta}$, using various special properties
of characters for the one part partition $(d)$. We
consider $\beta\vdash d$ with $l(\beta )=n$, and
let the number of parts of $\beta$ equal to $i$ be given by $N_i$, $i\geq 1$.
Thus $\sum_{i\geq 1}N_i =n$ and $\sum_{i\geq 1} iN_i = d$. We also let
$c_1=N_1-1$ and $c_i=N_i$, for $i\geq 2$, and
$$\B_{2j}=\sum_{i\geq 1}
i^{2j}c_i =-1+\sum_{i\geq 1}i^{2j}N_i = -1 + \sum \be_i^{2j},$$
 for $j\geq 1$,
a power sum for the partition, shifted by $1.$
Let $\xi_{2j}=\left[x^{2j} \right]\,\log( \sinh(x)/x)$,
and  let $\xi_\lambda=\xi_{\lambda_1}\xi_{\lambda_2} \dots$.
For any partition $\lambda =(\lambda _1,\ldots )$,
let $2\lambda =(2\lambda _1,\ldots )$.


\begin{theorem}\label{onepartg}\lremind{onepartg}Let 
$r=r^{\genus{g}}_{(d), \be}$.
For $g\geq 0$, and $\beta\vdash d$ with $n$ parts, \lremind{onepartg1, onepartg2}
\begin{eqnarray}\label{onepartg1}
H^{\genus{g}}_{(d),\beta } &=&
r! d^{r-1}\left[ t^{2g}\right]
\prod_{k \geq 1} \left(  \frac {\sinh(k t/2)} {k t/2} \right)^{c_k} \\
&=& \label{onepartg2} \frac{ r! d^{r-1} }{2^{2g} }
\sum_{\lambda\vdash g} \frac{\xi _{2\lambda}\B_{2\lambda}}{\left|
\Aut \lambda \right| }.
\end{eqnarray}
\end{theorem}

\noindent {\em Remarks.}
\newline \noindent
{\em 1.}
Equ.\ \eqref{onepartg1} is a generalization of a theorem of
Shapiro-Shapiro-Vainshtein: \cite[Thm.\ 6]{ssv} is the case $\be =
(1^d)$.  To our knowledge, \cite{ssv} contains the first appearance of
the generating series $(\sinh t)/t$ in connection with branched covers of curves (see
\cite[Sect.\ 3]{icm} for some of the subsequent connections, for
example through the Gopakumar-Vafa conjecture).  As noted there,
proofs of equivalent statements appear in \cite{jssv} and \cite{gssv},
but \cite{ssv} is the first interpretation in terms of Hurwitz
numbers.

\noindent
{\em 2.}
Equ.\ \eqref{onepartg1} also generalizes Theorem 3.2 of
\cite{gjcacti}:
\begin{equation*}
H^{\genus{0}}_{(d),\beta }=
 r! d^{r-1}.
\end{equation*}
(Note that, in ~\cite[Thm.\ ~3.2]{gjcacti},
the right hand side of the condition $t_1 +\cdots +t_m=n+1$ should
be replaced with $(m-1)n+1$.)

\proof We use the Hurwitz axioms (Prop.\ \ref{huraxioms}).
The group generated by any element of ${\mathcal{C}}_{(d)}$ acts
transitively on $\{ 1,\ldots ,d\}$.  But $H^{\genus{g}}_{\al,\beta }$ is a class function so,
for $\alpha =(d)$ and $r=r^{\genus{g}}_{(d),\be}=n-1+2g$, axiom~(H4) gives
\begin{equation*}
H^{\genus{g}}_{(d),\beta}=\frac 1 {\prod \be_j} 
\left[ {\mathsf{K}}_{\beta}\right]\left({\mathsf{K}}_{(2,1^{d-2})}\right) ^r
{\mathsf{K}}_{(d)}
= \frac 1 {d \prod \be_j}
\sum_{\lambda\vdash d}\eta(\lambda )^r\chi^{\lambda}_{(d)}
\chi^{\lambda}_{\beta}
\end{equation*}
from~(\ref{char1}) 
and ~(\ref{char3}), where
\begin{equation*}
\eta(\lambda )=\frac{\left|{\mathcal{C}}_{(2,1^{d-2})}\right| \chi^{\lambda}_{(2,1^{d-2})}}
{\chi^{\lambda}_{(1^d)}}=\sum_i{  \binom {\lambda _i}  2 }-\sum_i
{   \binom  {{\widetilde{\lambda}} _i}  2  },
\end{equation*}
and $\widetilde{\lambda}$ is the {conjugate} of $\lambda$. But
$\chi^{(d-k,1^k)}_{(d)} =(-1)^k,\; k=0,\ldots ,d-1$,
and $\chi^{\lambda}_{(d)}=0$,
for all other $\lambda$, and $\left|{\mathcal{C}}_{(d)} \right|  =(d-1)!$. Also
\begin{equation}\label{hookgf}
\sum_{k=0}^{d-1}\chi^{\left( d-k,1^k \right)}_{\beta}y^k=\prod_{i\geq 1} \left( 1-(-y)^i
\right)^{c_i},
\end{equation}
and $\eta((d-k,1^k))={\binom {d-k} 2}-{\binom {k+1} 2}=\binom d 2 -dk$. Thus
\begin{eqnarray*}
H^{\genus{g}}_{(d),\beta}&=&  \frac 1 {\prod \be_j} d^{r-1}
\sum_{k=0}^{d-1}\left(\frac{d-1}{2}-k\right) ^r(-1)^k
\chi^{(d-k,1^k)}_{\beta}\\
&=& \frac 1 {\prod \be_j} d^{r-1}\left[
\frac{t^r}{r!}\right]
\sum_{k=0}^{d-1}e^{\left(\frac{d-1}{2}-k\right) t} (-1)^k
\chi^{(d-k,1^k)}_{\beta}\\
&=& \frac 1 {\prod \be_j} r!d^{r-1}\left[ t^r\right]
e^{\frac{d-1}{2}t}\prod_{k \geq 1}\left( 1-e^{-kt}\right) ^{c_k},
\end{eqnarray*}
by substituting $y=-e^{-t}$ in~(\ref{hookgf}) above.
But $\sum_{h\geq 1}h c_h =d-1$ and $\sum_{h\geq 1} c_h =n-1$, so
we obtain
\begin{equation*}
H^{\genus{g}}_{(d),\beta}
= r! d^{r-1}\left[ t^{r-n+1}\right]
\prod_{k \geq 1} \left(  \frac {\sinh(kt/2)} {kt/2} \right)^{c_k}.
\end{equation*}
This yields \eqref{onepartg1}.  Applying the
logarithm,
$$\prod_{k\geq 1}
 \left(  \frac {\sinh(kx)} {kx} \right)^{c_k}
=\exp\left(\sum_{k\geq 1}c_k  \sum_{j\geq 1}
\xi _{2j}i^{2j}x^{2j}\right)
=\exp\left(\sum_{j\geq 1}\xi _{2j}\B_{2j}x^{2j}\right)=
\sum_{\lambda} \frac{\xi _{2\lambda}\B_{2\lambda}}{ \left|\Aut \lambda  \right| }
      x^{2\left|\lambda \right| } ,$$
where the sum is over all partitions $\lambda$.
Equ.\ \eqref{onepartg2} follows. \qed

Polynomiality is immediate from \eqref{onepartg2}; hence we have proved the
following.

\begin{corollary} \label{stronghomtrue} \lremind{stronghomtrue}The 
double Hurwitz numbers $H^{\genus{g}}_{(d), \be}$ satisfy
polynomiality and the Strong Piecewise Polynomiality Conjecture~\ref{stronghom}.
\end{corollary}
(The polynomials for $g \leq 5$ can be read off from Corollary~\ref{cphi}.)

Even more striking, the polynomial is divisible by $d^{n+2g-2}$, and
${H^\genus{g}_{(d),\be} }/{ r! d^{n+2g-2}}$
is a polynomial 
in the parts of $\be$ that is independent of the number of parts.
Examples may be seen in the next corollary, in which the one part double Hurwitz
numbers are expressed in terms of the $S_i$, which are polynomials in the parts
of $\beta$.

\begin{corollary}\label{cphi} \lremind{cphi}For 
$g\le 5$, explicit expressions for $H^{\genus{g}}_{(d),\beta}$
are given by 
\begin{eqnarray*}
H^{\genus{0}}_{(d),\be} &=&  (n-1)! d^{n-2} \\
H^{\genus{1}}_{(d), \be} &=& \frac{ (n+1)!}{24 }d^{ n} S_2 \\
H^{\genus{2}}_{(d), \be} &=& \frac{ (n+3)! d^{n+2}} { 5760  }
 \left(5S_2^2 - 2S_4\right) \\
H^{\genus{3}}_{(d),\be} &=& \frac { (n+5)! d^{n+4}} { 2^{10} \cdot 3^4 \cdot 5 \cdot 7  }
\left(16 S_6 -42 S_2 S_4 + 35 S_2^3 \right) \\
H^{\genus{4}}_{(d), \be} &=& \frac { (n+7)! d^{n+6}} {2^8}
\left(        -\frac{\B_8}{37800} + \frac{\B_2 \B_6}{17010} + \frac{\B_4^2}{64800}
      -\frac{\B_2^2 \B_4}{12960}+\frac{\B_2^4}{31104}
\right)
\\
H^{\genus{5}}_{(d), \be} &=& \frac { (n+9)! d^{n+8}} {2^{10}}
\left(  \frac{\B_{10}}{467775} - \frac{\B_2 \B_8}{226800} - \frac{\B_4 \B_6}{510300}
+\frac{\B_2^2 \B_6}{204120}+\frac{\B_2 \B_4^2}{388800}-\frac{\B_2^3 \B_4}{233280}
+\frac{\B_2^5}{933120}
\right)
\\
\end{eqnarray*}
\end{corollary}

Notice how constants associated to intersection theory on moduli spaces of 
low genus curves (such as $1/24$ and $1/5760$ for genus $1$ and $2$ respectively)
make their appearance.

In addition we note the following
attractive formula for the number of branched
covers of any genus and degree, with complete branching over two points.

\begin{corollary}\label{oneone} \lremind{oneone}
$$ H^\genus{g}_{(d),(d)} = (2g)! d^{2g-2} \left[ t^{2g} \right]
  \frac { \sinh (dt/2) }
{ \sinh (t/2) }
=d^{2g-2}\sum_{k=- \frac {d-1} 2}^{\frac {d-1} 2}k^{2g}.$$
\end{corollary}

\bpf The first equality is \eqref{onepartg1}, and the second comes
after straightforward manipulation. \epf

\subsection{From polynomiality to the symbol $\llangle \cdot \rrangle_g$, and
intersection theory on moduli} Theorem~\ref{onepartg}
strongly suggests the existence of an ELSV-type formula for one-part
double Hurwitz numbers, and even suggests the shape of such a formula.
In particular, we are in a much better position than we were for
single Hurwitz numbers when the ELSV formula \eqref{elsv} was
discovered.  At that point, polynomiality was conjectured
(\cite[Conj.\ 1.2]{gj2}, see also \cite[Conj.\ 1.4]{gjvn}).  Even
today, polynomiality has only been proved by means of the ELSV formula; no
character-theoretic or combinatorial reason is known.

In the one-part double Hurwitz case  we have much more.
\newline \noindent 
{\em (i)}
 We have a non-geometric proof of polynomiality.  
\newline \noindent 
{\em (ii)}  The polynomials
$P^{\genus{g}}_{1,n}(\be_1, \be_2, \dots)  = H^\genus{g}_{(d),\be}$
have an excellent description in terms of generating series. 
\newline \noindent 
 {\em (iii)} The polynomials are well-behaved as $n$ increases.  (More precisely, as described
before Corollary~\ref{cphi}, the polynomial $H^{\genus{g}}_{(d),\be}$ is
divisible by $d^{2g-2+n}$, and the
quotient $H^\genus{g}_{(d), \be} / r! d^{n+2g-2}$ is 
independent of $n$.)  Finally, 
\newline \noindent 
{\em (iv)} the polynomials may be seen, for non-geometric reasons, to satisfy the string and dilaton
equation.  (This will be made precise in Proposition~\ref{sd}.)

Hence we make the following geometric conjecture.  (The formula
\eqref{conj} is identical to \eqref{earlyconj} in the Introduction.)
The conjecture should be understood as: ``There exists a moduli space
$\oPic_{g,n}$ with the following properties...''. 

\begin{conjecture}[ELSV-type formula for one-part double Hurwitz numbers] \label{onepartconj} \lremind{onepartconj}For 
each $g \geq 0$, $n \geq 1$, $(g,n) \neq (0,1)$, $(0,2)$, 
 \lremind{conj}
\begin{equation}\label{conj}
H^{\genus{g}}_{(d), \be} = r^{\genus{g}}_{(d), \be}!  d 
\int_{\oPic_{g,n}}  \frac  { \Lambda_0 - \Lambda_2 + \cdots \pm \Lambda_{2g}}
{( 1 - \beta_1 \psi_1) \cdots (1 - \be_n \psi_n)},
\end{equation}
where $\oPic_{g,n}$, $\psi_i$, and $\Lambda_{2k}$ satisfy the following
properties.
\begin{itemize}
\item  {\em The space $\oPic$, and its fundamental class.}
  There is a
  moduli space $\oPic_{g,n}$, with a (possibly virtual) fundamental
  class $[ \oPic_{g,n} ]$ of dimension $4g-3+n$, and an open
  subset isomorphic to the Picard variety $\Pic_{g,n}$ 
  of the universal curve over
  $\cm_{g,n}$ (where the two fundamental classes agree).  
\item {\em Morphisms from $\oPic$.}  There is a forgetful
  morphism $\pi: \oPic_{g,n+1} \rightarrow \oPic_{g,n}$ (flat, of
  relative dimension 1), with $n$ sections $\si_i$ giving Cartier
  divisors $\De_{i,n+1}$ ($1 \leq i \leq n$).  Both morphisms behave
  well with respect to the fundamental class: $[\oPic_{g,n+1}] = \pi^*
  [\oPic_{g,n}]$, and $\De_{i,n+1} \cap \oPic_{g,n+1} \cong
  \oPic_{g,n}$ (with isomorphisms given by $\pi$ and $\si_i$),
  inducing $\De_{i,n+1} \cap [\oPic_{g,n+1}] \cong [\oPic_{g,n}]$.
\item {\em $\psi$-classes on $\oPic$.} There are $n$ line
  bundles, which are the cotangent spaces to the first $n$ points on
  $\Pic_{g,n}$, with first Chern classes $\psi_1$, \dots, $\psi_n$.  They
  satisfy $\psi_i = \pi^* \psi_i + \De_{i,n+1}$ ($i \leq n$) on
  $\oPic_{n+1}$ (the latter $\psi_i$ is on $\oPic_{n}$), and $\psi_i
  \cdot \De_{i,n+1} = 0$.
\item {\em $\Lambda$-classes.}  There are Chow (or cohomology)  classes
$\Lambda_{2k}$  ($k= 0, 1, \dots, g$) of codimension $2k$
on $\oPic_{g,n}$, which are pulled back from $\oPic_{g,1}$ (if $g>0$)
or $\oPic_{0,3}$; $\Lambda_0 = 1$.  The $\Lambda$-classes
are the Chern classes of a rank $2g$ vector bundle isomorphic to its
dual.  
\end{itemize}
  \end{conjecture}
  
  The suggestion that the $\Lambda$-classes are the Chern classes of a
  self-dual vector bundle is due to J. Bryan.  At the very least, one
  might expect that the $\Lambda_{2k}$ are tautological, given the
  philosophy that ``geometrically natural classes tend to be
  tautological'' (see {\em e.g.} \cite{notices}).

\noindent
{\em Remarks.}
\newline
\noindent
{\em 1.} Our 
motivation for this conjecture included {\em (a)}
the form of the ELSV-formula \eqref{elsv}, {\em (b)} 
the Piecewise Polynomiality Theorem~\ref{homthm},
{\em (c)} the presence of an $n$-pointed curve with line bundle (the
covering curve, with the points above $\infty$ and the pullback
of the line bundle ${\mathcal{O}}_{\C \proj^1}(1)$), and {\em (d)}
the remaining results of this section.  (In particular, the string
and dilaton equations, Prop.\ \ref{sd}, motivated the
conditions on the $\psi$-classes.)

\noindent
{\em 2.}  As pointed out in the introduction, the most speculative
part of this conjecture is the identification of the
$(4g-3+n)$-dimensional moduli space with a compactification of
$\Pic_{g,n}$; the evidence suggests a space of this dimension
with a morphism to $\cmbar_{g,n}$.  We suggest this space because of
{\em 1(c)} above, and the presence of the Picard variety in
symplectic field theory \cite{yasha}. 
 
\noindent
{\em 3.}  There are certainly other formulae for double Hurwitz
numbers not of this form, for example those involving integrals on the
space of relative stable maps.  However, to our knowledge, none of
these formulae explains polynomiality of one-part double Hurwitz
numbers, or the strong features of these polynomials.

\noindent
{\em 4.} A satisfactory proof would connect the geometry
of one-part double Hurwitz numbers with \eqref{conj}.

\noindent {\em 5.}  See Conjecture~\ref{bonusconj} relating
$\Lambda_{2g}$ to $\la_g$.

In analogy with Witten's notation \eqref{wittennotation},
we define\lremind{wittnot} $\llangle \tau_{b_1} \cdots \tau_{b_n} \Lambda_{2k} \rrangle_g $ by
\begin{equation}\label{wittnot}
\llangle \tau_{b_1} \cdots \tau_{b_n} \Lambda_{2k} \rrangle_g =
(-1)^{k} \left[ \be_1^{b_1} \cdots \be_n^{b_n} \right] \left(
\frac{  P^{\genus{g}}_{1,n}(\be_1, \dots, \be_n) }
{r^{\genus{g}}_{(d), \be}!d }\right) 
\end{equation}
if $(g,n) \neq (0,1)$, $(0,2)$, $\sum b_i + 2k = 4g-3+n$ and the $b_i$
are non-negative integers, and $\llangle \tau_{b_1} \cdots \tau_{b_n}
\Lambda_{2k} \rrangle_g := 0$ otherwise.  This definition makes sense
by Corollary~\ref{stronghomtrue}.  
Note that the symbol is symmetric in the $b_i$.
Conjecture~\ref{onepartconj}
then implies that \lremind{edwitten}
\begin{equation} \label{edwitten}
\llangle \tau_{b_1} \cdots \tau_{b_n} \Lambda_{2k} \rrangle_g
= \int_{\oPic_{g,n}} \psi_1^{b_1} \cdots \psi_n^{b_n} \Lambda_{2k}.
\end{equation}

\subsection{Generating series for $\llangle \cdot \rrangle_g$, and
the string and dilaton equations.}
This symbol has some remarkable properties which suggest geometric
meaning, in analogy with Witten's symbol $\langle \cdot \rangle_g$.  We
determine two expressions for a particular generating series for
this symbol, and then derive string and dilaton
equations, and prove Conjecture~\ref{onepartconj}
in genus $0$ and $1$.

Define $Q^{(i)}(x) = \sum_{j \geq 1} q_j j^i x^j$ for $i \geq 0$, so
$Q^{(0)}(x) = Q(x)$, defined just after \eqref{funceqw}, and
$$Q^{(i)}(x) = \left(  x \frac d {dx} \right)^i Q(x), \quad \quad i \geq 0.$$
Of course, we also have\lremind{diffrules}
\begin{equation}\label{diffrules}
x \frac d {dx}Q^{(i)}(x)=Q^{(i+1)}(x), \;\;\;\; i\geq 0.
\end{equation}
The first expression for the generating series follows directly
from the definition~(\ref{wittnot}) of the symbol $\llangle \cdot \rrangle_g$.

\begin{theorem}\label{Wittgfcn}
\lremind{Wittgfcn}For $g \geq 0$,  
\begin{equation*}
 x\frac{d}{dx}\sum_{n\geq 1} \frac{1}{n!}\sum_{k=0}^g (-1)^k
\sum_{b_1,\ldots ,b_n\geq 0}
\llangle \tau_{b_1} \cdots \tau_{b_n} \Lambda_{2k} \rrangle_g
\prod_{i=1}^n Q^{(b_i)}(x)
= \left. \HH_1^{\genus{g}}(x) \right|_{u=1}
\end{equation*}
\end{theorem}

\bpf
 From~(\ref{wittnot}), we have
\begin{eqnarray*}
LHS&=&x\frac{d}{dx}\sum_{n\geq 1} \frac{1}{n!}
\sum_{\beta_1,\cdots ,\beta_n\geq 0}
q_{\beta_1}\cdots q_{\beta_n}x^{\beta_1+\cdots +\beta_n}
\sum_{k=0}^g (-1)^k
\sum_{b_1,\dots ,b_n\geq 0} \beta_1^{b_1}\cdots \beta_n^{b_n}
\llangle \tau_{b_1} \cdots \tau_{b_n} \Lambda_{2k} \rrangle_g\\
&=&\sum_{d\geq 1}dx^d\sum_{n\geq 1}\frac{1}{n!}
\sum_{\beta_1+\cdots +\beta_n=d}q_{\beta_1}\cdots q_{\beta_n}
\sum_{k=0}^g (-1)^k
\sum_{b_1,\dots ,b_n\geq 0} \beta_1^{b_1}\cdots \beta_n^{b_n}
\llangle \tau_{b_1} \cdots \tau_{b_n} \Lambda_{2k} \rrangle_g\\
&=& \sum_{d\geq 1}dx^d\sum_{n\geq 1}
\sum_{\substack{{\beta\vdash d}\\{l(\beta) =n}}}
\frac{q_{\beta}}
{\left|\Aut\beta\right|}\sum_{k=0}^g (-1)^k
\sum_{b_1,\dots ,b_n\geq 0} \beta_1^{b_1}\cdots \beta_n^{b_n}
\llangle \tau_{b_1} \cdots \tau_{b_n} \Lambda_{2k} \rrangle_g\\
&=&RHS,
\end{eqnarray*}
giving the result.
\epf

The second expression for this generating series follows
from Theorem~\ref{onepartg}. To state this result requires some
more notation.
Define $v_{2j},\; f_{2j}$ by
$$\frac{2}{x}\sinh \left( \frac{x}{2} \right) 
=\sum_{j\geq 0}v_{2j}x^{2j},\;\;\;\;\;\;\;\;
  \frac{x}{2}\mathrm{cosech} \left( \frac{x}{2} \right)=
\sum_{j\geq 0}f_{2j}x^{2j}.$$
Then we have (see for example \cite[Sect.\ 1.41 and \ 9.6]{gr}),
$$v_{2j}=\frac{1}{2^{2j} (2j+1)!},\;\;\;\;\;\;\;\;
f_{2j}=\frac{1-2^{2j-1}}{2^{2j-1}(2j)!}B_{2j}, \;\;\;\; j\geq 0,$$
where $B_{2j}$ is a Bernoulli number ($B_0=1$, $B_2=1/6$, $B_4=-1/30$,
 $B_6=1/42$, $\dots$). As Bernoulli numbers alternate in sign
after $B_2$, note that $f_{2j}$ has sign $(-1)^j$, $j\geq 0$.
For a partition $\beta=(\beta_1,\ldots )$, let $Q^{(\beta)}(x)=\prod_{i\geq 1}
Q^{(\beta_i)}(x)$.

\begin{theorem}\label{Wittansx}\lremind{Wittansx}For $g\geq 0$,
$$\left. \HH_1^{\genus{g}}(x) \right|_{u=1}
=\sum_{k=0}^g f_{2k}\sum_{\substack{{\theta\vdash_0 g-k}
\\{l(\theta)\geq 1}}}
\frac{v_{2\theta}}{\left|\Aut\theta\right|}
\left(x\frac{d}{dx}\right)^{2g-2+l(\theta)}
Q^{(2\theta)}(x),$$
where $\theta\vdash_0 g-k$ means that $\theta$ is a partition
of $g-k$, with $0$-parts allowed.
\end{theorem}

\bpf
 From~(\ref{onepartg1}), we have
\begin{eqnarray*}
\left[ x^d \right] \left. \HH_1^{\genus{g}}(x) \right|_{u=1}
&=& \sum_{\beta \vdash d } \frac {H^{\genus{g}}_{(d), \beta}} {r! \Autbe} q_{\be} \\
&=& d^{2g-2}\left[ x^d t^{2g} \right]\frac{t/2}{\sinh(t/2)}
\exp\left(d\sum_{j\geq 1}q_j x^j\frac{\sinh(jt/2)}{jt/2}\right)\\
&=& d^{2g-2} \sum_{k=0}^g f_{2k}  \left[ x^d t^{2g-2k} \right]\exp
  \left(  d \sum_{i \geq 0} v_{2i} t^{2i} Q^{(2i)}(x) \right) \\
&=& d^{2g-2} \sum_{k=0}^g f_{2k} \left[ x^d \right]\sum_{\theta\vdash_0 g-k}
d^{l(\theta)}\frac{v_{2\theta}}{\left|\Aut\theta\right|}
Q^{(2\theta)}(x), 
\end{eqnarray*}
and the result follows. 
\epf

By comparing the two generating series expressions given 
in Theorems~\ref{Wittgfcn}
and~\ref{Wittansx}, we obtain an explicit expression
for $\llangle\cdot\rrangle_g$, in the following result.

\begin{corollary}\label{Wittcor}
For all $b_1$, \dots, $b_n$, $k$, $g$, with $b=(b_1,\ldots ,b_n)$,
$$\llangle \tau_{b_1} \cdots \tau_{b_n} \Lambda_{2k} \rrangle_g
=\left|\Aut b \right|f_{2k}(-1)^k\sum_{\substack{{\theta\vdash_0 g-k}\\
{l(\theta)=n}}}
\frac{v_{2\theta}}{\left|\Aut\theta\right|}\left[ Q^{(b_1)}\cdots Q^{(b_n)}
\right]
\left(x\frac{d}{dx}\right) ^{2g-3+n} Q^{(2\theta)}(x).$$
\end{corollary}

\bpf
Compare Theorems~\ref{Wittgfcn} and~\ref{Wittansx}.
Now $\llangle \tau_{b_1} \cdots \tau_{b_n} \Lambda_{2k} \rrangle_g$ is
symmetric in the $b_i$'s, so each of the $n!/\left|\Aut b \right|$ distinct
reorderings of $b$ contributes equally to the coefficient of the
monomial $Q^{(b_1)}\ldots Q^{(b_n)}$. 
The result
follows from~(\ref{diffrules}).
\epf

{}From~Corollary~\ref{Wittcor}
and applying~(\ref{diffrules}), we can obtain a
great deal of information about values of $\llangle\cdot\rrangle_g.$
For example, we immediately have the following non-negativity result.
(The analogous result for the Witten symbol $\langle \tau_{a_1} \cdots
\tau_{a_n} \la_k \rangle_g \geq 0$ is non-obvious.)

\begin{corollary}[Non-negativity]
For all $b_1$, \dots, $b_n$, $k$, $g$, $\;\;\;\;
\llangle \tau_{b_1} \cdots \tau_{b_n} \Lambda_{2k} \rrangle_g \geq 0$.
\end{corollary}

We can also prove that the symbol $\llangle \cdot \rrangle_g$ satisfies
the string and dilaton equations, in the following result.

\begin{proposition}[String and dilaton equations] \label{sd} \lremind{sd}
{} \quad
\begin{enumerate} 
\item[(a)] {\em (string equation)}  The following equation holds, except
when $g=k=1, n=0$.
$$
\llangle \tau_0 \tau_{b_1} \cdots \tau_{b_n} \Lambda_{2k} \rrangle_g
= \sum_{i=1}^n \llangle \tau_{b_1} \cdots \tau_{b_{i-1}} \tau_{b_i-1}
\tau_{b_{i+1}} \cdots
\tau_{b_n} \Lambda_{2k} \rrangle_g$$
In the exceptional case, we have $\llangle \tau_0 \Lambda_{2}\rrangle_1=1/24$.
\item[(b)] {\em (dilaton equation)} The following equation holds.
  $$
  \llangle \tau_1 \tau_{b_1} \cdots \tau_{b_n} \Lambda_{2k}
  \rrangle_g = (2g-2+n) \llangle \tau_{b_1} \cdots
  \tau_{b_n} \Lambda_{2k} \rrangle_g$$
\end{enumerate}
\end{proposition}

Note that, by the usual proofs of the string and dilaton equation
(see for example \cite[Sect.\ 1]{looijenga} or \cite[Ch.\ 25]{MirSym}),
this proposition would be implied by Conjecture~\ref{onepartconj}.

\bpf
For formal power series in the variables $Q^{(i)}:=Q^{(i)}(x)$, $i\geq 0$,
we define the partial differential operators\lremind{september}
\begin{equation}
\label{september}\Delta_{-1}=\sum_{i\geq 0 } Q^{(i+1)}\frac{\partial}{\partial Q^{(i)}}
,\;\;\;\; \Delta_0=\sum_{i\geq 0} Q^{(i)}\frac{\partial}{\partial Q^{(i)}}.
\end{equation}
Note that we have the operator identity $\Delta_{-1}=x\frac{d}{dx}$, as well as
\begin{equation}\label{opone}
\frac{\partial}{\partial Q^{(i)}}\Delta_{-1}=\Delta_{-1}
\frac{\partial}{\partial Q^{(i)}}+\frac{\partial}{\partial Q^{(i-1)}},
\end{equation}
and
\begin{equation}\label{optwo}
 Q^{(i)}\Delta_{-1}=\Delta_{-1}
 Q^{(i)}- Q^{(i+1)}.
\end{equation}
Now, multiplying~(\ref{opone}) on the left by $Q^{(i)}$, applying~(\ref{optwo}),
and then summing over $i\geq 0$, we obtain
\begin{equation}\label{opthree}
\Delta_0 \Delta_{-1}= \Delta_{-1} \Delta_0.
\end{equation}
Also, applying~(\ref{opone}) repeatedly with $i=0,1$, we obtain
\begin{equation}\label{opfour}
\frac{\partial}{\partial Q^{(0)}}\Delta_{-1}^m=\Delta_{-1}^m \frac{\partial}{\partial Q^{(0)}},
\;\;\;\;\;\;\;\;\frac{\partial}{\partial Q^{(1)}}\Delta_{-1}^m=\Delta_{-1}^m \frac{\partial}{\partial Q^{(1)}}
+m\Delta_{-1}^{m-1} \frac{\partial}{\partial Q^{(0)}}.
\end{equation}
Now let $x\frac{d}{dx}\Psi_{\genus{g}}(x)=\left. \HH_1^{\genus{g}}(x) \right|_{u=1}$.
Then, from Theorem~\ref{Wittansx} and~(\ref{opfour}), we obtain
$$\left( \frac{\partial}{\partial Q^{(0)}}-\Delta_{-1}\right)\left(
\Psi_{\genus{g}}+\delta_{g,1}\frac{Q^{(0)}}{24}\right) =0,$$
and thus the string equation holds with the
given exceptional value, from Corollary~\ref{Wittcor}.

Also, from Theorem~\ref{Wittansx} and~(\ref{opthree}),~(\ref{opfour}), we obtain
$$\left( \frac{\partial}{\partial Q^{(1)}}-\Delta_0-(2g-2)\right)\Psi_{\genus{g}}=0,$$
and thus the dilaton equation holds, from Corollary~\ref{Wittcor}.
\epf

\subsubsection{Virasoro constraints?}
In the case of the moduli space of curves, the string and dilaton
equations are essentially the first two Virasoro constraints (see for
example \cite[Sect.\ 25.2]{MirSym}).  It is natural then to ask
whether there is a full set of Virasoro constraints.  Even in the case
of single Hurwitz numbers, this is not known.  However, in the single
Hurwitz number case, the highest-degree terms (of the polynomial
defined in \eqref{eHCP}) are polynomials with coefficients of the form
$\langle \tau_{a_1} \cdots \tau_{a_m} \rangle_g$ ({\em i.e.} with no
$\la$-class), which {\em do} satisfy Virasoro constraints, by Witten's
conjecture (Kontsevich's theorem)~\cite{witten,kont}.  (Indeed, this
idea led to Okounkov and Pandharipande's proof of Witten's conjecture
\cite{op1}.)  Thus one may ask a weaker question: are there Virasoro
constraints on the asymptotics of one-part double Hurwitz numbers,
{\em i.e.} on $\llangle \tau_{b_1} \cdots \tau_{b_n} \rrangle_g
:= \llangle \tau_{b_1} \cdots \tau_{b_n} \Lambda_0 \rrangle_g$?
Given Conjecture~\ref{onepartconj}, this is the analogue of Witten's
conjecture on the compactified Picard variety.  

We have not yet been able to produce a set of Virasoro constraints, but
our partial results suggest additional hidden structure, so we report them here
without proof.

For formal power series in the variables $Q^{(i)}:=Q^{(i)}(x)$, $i\geq 0$,
we define the partial differential operators
$\Delta_{-1}$, $\Delta_0$ (as in \eqref{september}),
$$\Delta_0^{\prime}=\sum_{i\geq 0} i Q^{(i)}\frac{\partial}{\partial Q^{(i)}},
\;\;\;\; 
\Delta^{\prime}_{1}=\sum_{i\geq 1 } i Q^{(i-1)}\frac{\partial}{\partial Q^{(i)}}
,\;\;\;\; \Delta_{1}^{\prime \prime}=\sum_{i\geq 1 } i^2 Q^{(i-1)}\frac{\partial}{\partial Q^{(i)}}.$$
Now let 
$$\Psi=\sum_{g\geq 0}\left.\Psi_g\right| _{k=0},$$
where  $\Psi_g$ is as above. (Note
that the value of $g$ is recoverable from the partition condition
on the monomials.)
Then the string equation translates to an 
annihilator $A_{-1}$ for $\Psi$ (up to
initial conditions), where
$$A_{-1} =
\frac{\partial}{\partial Q^{(0)}}-\Delta_{-1}.$$
The dilaton equation translates to an 
annihilator $A_0$ for $\Psi$, where
$$A_0=
\frac{\partial}{\partial Q^{(1)}}-\frac{1}{2}\left(\Delta_0+
  \Delta_0^{\prime}-1\right).$$
It is an easy computation that
$$\left[A_{-1},A_0\right]=\frac{1}{2}A_{-1}.$$
Now \cite[p. 5689]{iz} suggests letting 
$B_{-1}=-\frac{1}{2}A_{-1}$ and 
$B_0=-2A_0$, so that
$$\left[B_{-1},B_0\right]=-B_{-1}.$$
We then sought a candidate $B_1$ (analogous to Witten's $L_1$)
involving a term of the form 
$\frac{\partial}{\partial Q^{(2)}}$.  There
are many such annihilators, and the simplest 
we found was
$$A_{1}=4\frac{\partial}{\partial Q^{(2)}}+\frac{1}{3}
\frac{\partial^3}{\partial {Q^{(0)}}^3}+2\left(\Delta_1^{\prime\prime}+
\Delta_1^{\prime}\right)\frac{\partial}{\partial Q^{(1)}}
-12\frac{\partial}{\partial Q^{(2)}}\frac{\partial}{\partial Q^{(1)}}.$$
However, we have been unable to find a candidate $B_1$ satisfying
the desired Virasoro commutation relations with $B_0$ and $B_{-1}$.

\subsubsection{Verifying Conjecture~\ref{onepartconj}  in low genus.}

\begin{proposition}
  Conjecture~\ref{onepartconj} is true in genus
  $0$, taking $\oPic_{0,n} = \cmbar_{0,n}$, and in genus $1$, taking
  $\oPic_{1,n} = \cmbar_{1,n+1}$ and $\Lambda_2 = \pi^* [pt]/24$,
  where $pt$ is the class of a point on $\oPic_{1,1}$, and $\pi$ is
  the morphism $\oPic_{1,n} \rightarrow \oPic_{1,1}$.
\end{proposition}

We have two proofs, neither of which is fully satisfactory (in the
sense of Remark~4 after Conjecture~\ref{onepartconj}).  First, the geometric arguments of \cite{v}
apply with essentially no change; this argument is omitted for the
sake of brevity.  The following second proof is purely combinatorial.

\bpf For genus $0$, if $n \geq 3$, then
$$ r^{\genus{0}}_{(d), \be}!  d 
\int_{\cmbar_{0,n}}  \frac   1
{( 1 - \beta_1 \psi_1) \cdots (1 - \be_n \psi_n)}
= r^{\genus{0}}_{(d), \be}!  d
( \be_1 + \cdots + \be_n)^{n-3} $$
(using the string equation; see for example \cite[Ex.\ 25.2.8]{MirSym}),
so we are done by Corollary~\ref{cphi}.

For genus $1$, we will prove \eqref{edwitten}.  We need only prove
the base cases
$\llangle \tau_2 \Lambda_0 \rrangle_1 = \frac 1 {24}$ and $\llangle
\tau_0 \Lambda_2 \rrangle_1 = \frac 1 {24}$ (obtained by unwinding
Corollary~\ref{cphi}), as the rest follow by the string and dilaton
equation.  The first is
$$
\int_{\cmbar_{1,2}} \psi_1^2 = \frac 1 {24},
$$ which is well-known ({\em e.g.}
\cite[Exer.\ 25.2.9]{MirSym}; combinatorialists may prefer to extract
it from the ELSV formula \eqref{elsv}), and the second is immediate from the
definition of $\Lambda_2$. \qed

\subsection{Explicit formulae for $\llangle \cdot \rrangle_g$}
We can also determine explicit formulae for many instances of
this symbol, just as such formulae have been given for 
Witten's symbol $\langle \cdot \rangle_g$, most
notably by Faber and Pandharipande.  In particular:

\noindent {\em Integrals over $\cmbar_{g,1}$.}
 It is a straightforward consequence of Witten's conjecture that
\lremind{babywitten}
\begin{equation}
\label{babywitten}
\langle \tau_{3g-2} \rangle_g = \frac 1 {24^g g!}
\end{equation}
(see for example just before (4) in \cite{nonvanishing}).
Also, \cite[equ.\ (5)]{hodgegw}: \lremind{wittenesque}
\begin{equation} \label{wittenesque}
\int_{\cmbar_{g,1}} \psi_1^{2g-2} \la_g = \frac{ 2^{2g-1}-1 } {2^{2g-1} (2g)!} \left| B_{2g} \right|.
\end{equation}
Generalizing both of these statements is
\cite[Thm.\ 2]{hodgegw}: \lremind{fp1}
\begin{equation}
\label{fp1}
1+\sum_{g\geq 1}\sum_{i=0}^g t^{2g} k^i
\int_{\cmbar_{g,1}} \psi_1^{2g-2+i} \la_{g-i} = \left(\frac{t/2}{\sin(t/2)} \right)^{k+1}.
\end{equation}

\noindent
{\em The ``$\la_g$-theorem''.}  The main theorem of \cite{lambdag} is:
\lremind{lambdagthm}
\begin{equation} \label{lambdagthm}
\langle \tau_{b_1} \cdots \tau_{b_n} \la_g \rangle = 
\int_{\cmbar_{g,n}} \psi_1^{b_1} \cdots \psi_n^{b_n} \la_g = \binom
{2g+n-3}{b_1, \dots, b_n} \mathbf{b}_g. \end{equation}
(This was
first conjectured in \cite[equ.\ (16)]{getzlerpandharipande}.  More precisely, it
was shown to be a consequence of the Virasoro conjecture for constant
maps to $\C \proj^1$.)  The constant $\mathbf{b}_g$ can be evaluated
using
\eqref{wittenesque} by taking $(b_1, b_2, b_3, \dots) = (3g-2, 0, 0, \dots)$.

We now deduce analogues and generalizations of these results for double Hurwitz
numbers.  A proof of Conjecture~\ref{onepartconj} would
thus give these results important geometric meaning.

In analogy with the $\la_g$-theorem \eqref{lambdagthm}, we have the
following result, which follows immediately from Theorem~\ref{onepartg}
(in the same way as did Cor.\ \ref{cphi}).
\begin{proposition}
\begin{equation}
\label{brunch}
\llangle \tau_{b_1} \cdots \tau_{b_n} \Lambda_{2g} \rrangle_g
= \left[  \beta_1^{b_1} \cdots \beta_n^{b_n} \right] \left( {\mathbf{c}}_g d^{r-2} + 
   \mbox{higher terms in $\be$'s} \right)\end{equation}
\lremind{brunch}($d = \sum \be_j$), where ${\mathbf{c}}_g$ depends only on $g$.  As $b_1 + \cdots + b_n = 2g-3+n = r-2$, we have
$$
\llangle \tau_{b_1} \cdots \tau_{b_n} \Lambda_{2g} \rrangle_g
= \binom {2g-3+n} {b_1, \dots, b_n} {\mathbf{c}}_g$$
for some constant $\mathbf{c}_g$.
\end{proposition}

We note that \eqref{brunch} is analogous to the version
\cite[equ.\ (18)]{lambdag} of the $\la_g$-theorem used in the proof of
Faber and Pandharipande.

In analogy with \eqref{babywitten}, we have\lremind{plant1}
\begin{equation} \label{plant1}
\llangle \tau_{4g-2} \rrangle_g = \frac 1 {2^{2g} (2g+1)!}.
\end{equation}
In analogy with \eqref{wittenesque}, we have\lremind{plant2}
\begin{equation} \label{plant2}
\llangle \tau_{2g-2} \Lambda_{2g} \rrangle_g =
\frac { (-1)^g (1 - 2^{2g-1})} { 2^{2g-1} (2g)!} B_{2g} 
=  \frac { 2^{2g-1} - 1} {2^{2g-1}(2g)!} \left| B_{2g} \right|.
\end{equation}
Thus we have evaluated ${\mathbf{c}}_g$ in the
previous Proposition.  Remarkably, it is the same constant appearing
in Faber and Pandharipande's expression \eqref{wittenesque}, 
leading us to speculate the following.

\begin{conjecture} \label{bonusconj}\lremind{bonusconj}There 
is a structure morphism $\pi:  \oPic_{g,n}
\rightarrow \cmbar_{g,n}$,
and  $\pi_* \Lambda_{2g} = \la_g$.
\end{conjecture}

Generalizing \eqref{plant1} and \eqref{plant2}, we have the following
result, in analogy with (but not identical to) \eqref{fp1}.

\begin{proposition} For $g \geq 1$, and $k=0$, \dots, $g$,
$$  \llangle \tau_{b_1} \Lambda_{2k} \rrangle _g =
(-1)^k f_{2k} v_{2g-2k} =
    \frac{(-1)^k f_{2k}}{2^{b_1-2g+2}(b_1-2g+3)!}$$
for $b_1+2k = 4g-2$.
Equivalently,
$$
1+\sum_{g\geq 1} t^{2g} \sum_{k=0}^g x^{2k}
\llangle \tau_{b_1} \Lambda_{2k} \rrangle_g
= \frac{ x \sinh (t/2) }{ \sin (xt/2) }.$$
\end{proposition}

\bpf
This follows immediately from Corollary~\ref{Wittcor}, since the only
choices of $\theta$ in the summation are partitions with a single part.
\epf

This result can be extended to expressions for terms with more $\tau$'s.   
For example, part (a) of the following proposition gives
a closed-form expression for any term involving two $\tau$'s.
There are also formulae for any mixture of $\tau_2$'s and $\tau_3$'s 
(where the number of $\tau_3$'s is held fixed);  
the first three examples are parts (b)--(d) below.
We know of no analogue for $\langle \cdot \rangle_g$.

\begin{proposition} \quad \newline
\begin{enumerate}
\item[(a)]  For $k=0$, \dots, $g$, and $g \geq 2$,
  $$\llangle \tau_{b_1} \tau_{b_2} \Lambda_{2k} \rrangle _g =
    \frac{(-1)^k f_{2k}}{2^{2g-2k+1}(2g-2k+2)!} \sum_{\text{$i>0$ odd}}
    \binom {2g-2k+2} i \left( 
 \binom {2g-1}{ b_1+1-i} + \binom {2g-1 }{b_2+1-i} \right)$$
for $b_1 + b_2 = 4g-2k-1$.   
\item[(b)]
For $k=0$, \dots, $g$, and $g\geq 1$, except $(k,g)=(1,1)$,
$$ \llangle \tau_2^{4g-3-2k} \Lambda_{2k} \rrangle _g = \frac{(-1)^k
  f_{2k}}{24^{g-k}(g-k)!} (6g-7-2k)!!$$
where $(2m-1)!! =
(2m-1)(2m-3) \cdots (3)(1)$ for $m$ a positive integer, and $(-1)!!=1$.
\item[(c)]
For $k=0$, \dots, $g$, and $g\geq 2$, except $(k,g)=(2,2)$,
$$    \llangle \tau_2^{4g-5-2k} \tau_3 \Lambda_{2k} \rrangle _g =
    \frac{(-1)^k f_{2k}}{24^{g-k}(g-k)!} (6g-7-2k)!! \frac{6g-4-4k}{3}.$$
\item[(d)]
For $k=0$, \dots, $g$, and $g\geq 2$, except $(k,g)=(1,2),(2,2),(3,3)$,
\begin{eqnarray*}
\lefteqn{    \llangle \tau_2^{4g-7-2k} \tau_3^2 \Lambda_{2k} \rrangle _g =}
\\
& &    \frac{(-1)^k f_{2k}}{24^{g-k}(g-k)!} (6g-9-2k)!!
    \frac{(3g-4-k)((6g-4-4k)(6g-7-4k)-(6g-2-6k))}{9}\end{eqnarray*}
\end{enumerate}
\end{proposition}

\bpf
These results all follow from Corollary~\ref{Wittcor} in a routine way,
using Leibniz's Rule.
For part (a), the only choices of $\theta$ in the summation are partitions with
two parts. For parts (b)--(d), all parts of $\theta$ must be $0$'s or $1$'s
only.
\epf

\subsection{A genus expansion ansatz for $\llangle \cdot \rrangle_g$ in the
  style of Itzykson and Zuber}
\label{new} \lremind{new}We next prove an analogue of the genus
expansion ansatz of Itzykson and Zuber for intersection numbers on the
moduli space of curves \cite[(5.32)]{iz}.  The
Itzykson-Zuber Ansatz was proved in \cite{eyy} and later in
\cite[Thm.\ 3.1]{gjvk}; the latter proof (and generalization) is similar
in approach to the argument in this paper.

\begin{theorem}[Genus Expansion Ansatz] \label{ResultI}
\lremind{ResultI}For $g \geq 0$,  \lremind{ResultIeq}
\begin{equation}
\label{ResultIeq}
\left. \HH_1^{\genus{g}}(x) \right|_{u=1} =
\sum_{k=0}^g f_{2k} \sum_{\theta \vdash g-k} \frac {v_{2 \theta}} { \left|
\Aut \theta \right|} \left(  x \frac d {dx} \right)^{2g-2+l(\theta)} \left( 
\frac {Q^{(2 \theta)} (w)} { 1 - Q^{(1)}(w) } - \delta_{\theta \emptyset} 
\right).
\end{equation}
\end{theorem}

\noindent {\em Remarks.} 

\noindent
{\em 1.}  Unlike the Itzykson-Zuber ansatz, this result
has explicitly computable coefficients.

\noindent
{\em 2.} $u=1$ is convenient here.  

\noindent
{\em 3.}  $w$ and $x$ are related by
\eqref{funceqw}.


\bpf
{}From Theorem~\ref{Wittansx}, we obtain
$$\left[ x^d \right]\left. \HH_1^{\genus{g}}(x) \right|_{u=1}
=\sum_{k=0}^g f_{2k}\sum_{\theta\vdash g-k}
\frac{v_{2\theta}}{\left|\Aut\theta\right|}
d^{2g-2+l(\theta)}\left[ x^d \right]
Q^{(2\theta)}(x)\sum_{m\geq 0}\frac{d^mQ^{(0)}(x)^m}{m!},$$
where $m$ is the number of $0$'s in the partition with $0$-parts allowed.
The result follows from
Lagrange's Implicit Function Theorem~\ref{lagthm}\eqref{lagthm2}
(with $u=1$ from \eqref{lagpartial}, \eqref{opident}).
\epf

To obtain explicit results from Theorem~\ref{ResultI}, we 
modify~(\ref{diffrules}), using~(\ref{lagpartial}) and~(\ref{opident}), to obtain
\begin{equation*}
x\frac{d}{dx}Q^{(i)}(w)=\frac{Q^{(i+1)}(w)}{1-Q^{(1)}(w)},\;\;\;\; i\geq 0.
\end{equation*}
For example, with $g=0$ in Theorem~\ref{ResultI}, and $i=0$ above, we obtain
\begin{equation}\label{ans0}
x\frac{d}{dx}\left. \HH_1^{\genus{0}}(x) \right|_{u=1}= Q^{(0)}(w)=Q(w).
\end{equation}
With $g=1$ and $i=1$, we obtain
\begin{equation}\label{ans1}
\left. \HH_1^{\genus{1}}(x) \right|_{u=1}= \frac{1}{24}\left(
Q^{(3)}(w)\mu(w)+Q^{(2)}(w)^2\mu(w)^2-\mu(w)+1\right),
\end{equation}
since, with $u=1$, we have $\mu(w)=1/(1-Q^{(1)}(w))$.

\noindent
{\em Remarks.}

\noindent {\em 1.}
 In genus $0$, there is a direct connection between the
generating series for single
Hurwitz numbers (with a partition $\be$) and one-part double Hurwitz
numbers.  More precisely, these generating series
are identical, under $\frac {j^j} {j!} p_j \leftrightarrow
q_j$ and $s \leftrightarrow w$. Here $s$ is the solution to the
functional equation
\begin{equation*}
s=x\, e^{\phi_0(s)},
\end{equation*}
and $\phi_i(x) = \sum_{j\geq 1} \frac{j^{j+i}}{j!}p_j x^j$,
as described
in \cite[Sect.\ 2.3]{gjvk}, so, for example, $\phi_0(x)
\leftrightarrow Q^{(0)}(x)=Q(x)$.  This is a purely formal
statement that the formula for single Hurwitz numbers and that for
one-part double Hurwitz numbers are ``essentially'' the same in genus
$0$.  We do not know if there is any geometric or combinatorial reason for this
coincidence.

\noindent {\em 2.} More generally, in arbitrary genus, there is also such
a connection. In this case, the direct
analogue of Theorem~\ref{Wittgfcn} holds for the single Hurwitz
number series, under $\phi_i(x)\leftrightarrow Q^{(i)}(x)$ and
 $\langle \cdot \rangle_g\leftrightarrow\llangle \cdot \rrangle_g$,
as described in \cite[Sect.\ 2.4]{gjvk}. However, there is no
analogue of Theorem~\ref{Wittansx} that we know of for the
single Hurwitz number series. From this point of view, the
Itzykson-Zuber Ansatz for the single Hurwitz number series
is the analogue of the form given by Theorem~\ref{ResultI}
under $s\leftrightarrow w$.

\noindent {\em 3.} We note that the substitution for $x$ by a series in $w$
specified by the 
functional equation~(\ref{funceqw}) 
is the key technical device used in Section~\ref{marb},
as considered in~(\ref{changevars}). However, the
approach in Section~\ref{marb} is completely different from that of the
present section, so the appearance of $w$ again suggests that
it is significant, and that a geometric or combinatorial
explanation for this would be interesting.

\noindent
{\em Caution.}  The results of this section, especially the Genus
Expansion Ansatz and the string and dilaton equations, seem to lead
inescapably to Conjecture~\ref{onepartconj}, but this is not quite the
case.  The simple structure of the polynomials $P^{\genus{g}}_{1,n}$
allows other possible statements as well.  For example, the correct
statement might be:
$$
H^{\genus{g}}_{(d), \be} =
 r^{\genus{g}}_{(d), \be}! \int_{\oPic'_{g,n+1}}
\frac { \Lambda'_0 - \Lambda'_2 + \cdots \pm \Lambda'_{2g}} {( 1 -
  \beta_1 \psi_1) \cdots (1 - \be_n \psi_n)},$$
 where the space $\oPic'_{g,n+1}$
and classes $\psi'_1$, \dots, $\psi'_n$, $\Lambda'_0$, \dots, $\Lambda'_{2g}$
satisfy the itemized hypotheses of Conjecture~\ref{onepartconj}.
(Notice that there is no ``$d$'' in the numerator, as there is
in Conj.\ \ref{onepartconj}.)
We use primes to indicate
that these objects need not be the same as in Conjecture~\ref{onepartconj}.

The $(n+1)$-st point should correspond to $\al$ (the point mapping to
$0$ in the target $\C \proj^1$).  $\oPic'_{g,n+1}$ should admit an action of $\mathfrak{S}_n$
(permuting the points corresponding to $\be$), but not necessarily $\mathfrak{S}_{n+1}$.  (This
insight comes from M. Shapiro, who has suggested that the correct
moduli space of curves for the double Hurwitz problem in general
should have two ``colors'' of points, one corresponding to $\al$,
and one corresponding to $\be$.) The string and dilaton equations are again satisfied.

\section{A symmetric function description of the Hurwitz generating series}
\label{character2} \lremind{character2}In
this section, we use character theory again,  to
give a good description of the double Hurwitz generating series $H.$
This gives
algebraic, rather than geometric, insight into Hurwitz numbers,
and thereby give a means of producing explicit formulae, for
example extending results of Kuleshov and M. Shapiro~\cite{ks}.

For the purposes of this section, we regard each of the indeterminates $p_k$ and
$q_k$ as power sum symmetric functions in two sets of indeterminates, one for
$p_k$ and the other for $q_k.$ This may be done since the power sum symmetric
functions in an infinite set of indeterminates are algebraically
independent. The following result gives such an expression, stated in terms of
symmetric functions. Let $s_{\lambda}(p_1,p_2,\ldots )$ be the Schur symmetric
function, written as a polynomial in the power sum symmetric functions
$p_1,p_2,\ldots$.
This is the generating series for the irreducible $\Sg_d$-characters
$\left(\chi^\alpha_\lambda\colon\lambda\vdash d\right)$ with respect to
the power sum symmetric functions.  For $\alpha\vdash d$, the expression, and
its inverse, is
\begin{equation}\label{spps}
s_{\alpha} =\frac{1}{d!}\sum_{\lambda\vdash d} \left|{\mathcal{C}}_{\lambda}\right|
\chi^{\alpha}_{\lambda}p_{\lambda},\;\;\;\;
p_{\alpha}=\sum_{\lambda\vdash d} \chi^{\lambda}_{\alpha}s_{\lambda}.
\end{equation}
{}From the expression for $H$ that we give next, we shall determine how
$H^{\genus{g}}_{\al,\be}$ depends on $g$ for fixed $\al$ and $\be.$

\begin{theorem}\label{Schurseries} \lremind{Schurseries}Let
\begin{equation*}
Z =1+ \sum_{d\geq 1}z^d\sum_{\lambda\vdash d}
e^{\eta (\lambda )t}s_{\lambda}(p_1t^{-1},p_2t^{-1},\ldots )
s_{\lambda}(q_1 u\, t^{-1},q_2 u\, t^{-1},\ldots ) .
\end{equation*}
Then
$\left.H\right| _{y=t^2} =t^2 \log Z$
and $\left. \tH \right|_{y=t^2} = t^2 Z$.
\end{theorem}

\proof Following the method of proof of Theorem~\ref{onepartg}, we have
\begin{equation}\label{discfactn}
\frac{\tH^{\,
(g)}_{\alpha,\beta}}{\Autalbe}
=\frac{\left|{\mathcal C}_{\beta}\right|}{d!}
\left[ {\mathsf{K}}_{\beta}\right]\left({\mathsf{K}}_{(2,1^{d-2})}\right) ^r
{\mathsf{K}}_{\alpha}
= \frac{\left|{\mathcal C}_{\al}\right|\cdot\left|{\mathcal C}_{\beta}\right|}{d!^2}
\sum_{\lambda\vdash d}\eta(\lambda )^r\chi^{\lambda}_{\alpha}
\chi^{\lambda}_{\beta}.
\end{equation}
Now multiply by $p_\alpha q_\beta u^{l(\beta )} z^d t^r/r!$, and sum over
  $\alpha,\beta\vdash d$,  $d\ge0$, and $r\geq 0$ (this number
is $0$ unless $r$ has the same parity as $l(\alpha )+l(\beta )$),
using~(\ref{spps}), to obtain
$$
\sum_{r\geq 0}
\sum_{d\geq 0}z^d\sum_{\alpha ,\beta\vdash d}
p_{\alpha}q_{\beta} u^{l(\beta )} t^r
\frac{\tH^{\, (g)}_{\alpha,\beta}} { \Autalbe r!}
=1+\sum_{d\ge1}z^d \sum_{\lambda\vdash d} e^{\eta(\lambda)t}
s_\la(p_1,\dots)s_\la(q_1 u,\ldots).
$$
This series is an exponential generating series in both $z$, marking
sheets, and $t$, marking transposition factors (we have divided by
both $r!$ and $d!$, the latter in the Hurwitz axioms Prop.\ \ref{huraxioms}).
To transform the exponent of $t$ from number of transposition factors
to genus, we apply the
substitutions $p_i\mapsto p_i t^{-1}, q_i \mapsto q_i t^{-1}$, $i\geq
1$, to obtain
$p_\alpha q_\beta t^r \mapsto
p_\alpha q_\beta t^{r-l(\alpha)-l(\beta)}
=
p_\alpha q_\beta t^{2g-2}$,
from~(\ref{rh2}),
and the result now follows.
(Note that $Z$ is clearly an even series in $t$ since
$\eta(\widetilde{\lambda)}=-\eta(\lambda)$ and
$s_{\tilde{\lambda}}(p)=s_{\lambda}(-p).$)
\qed


\subsection{Expressions for $H^{\genus{g}}_{\al,\be}$ for varying $g$
and fixed $\al,\be$}

Theorem~\ref{Schurseries} may be used to obtain $H_{\al,\be}^{\genus{g}}$
for fixed $\al$, $\be$.  The expressions are linear combinations
of $g$th powers of non-negative integers. In particular, the results
of Kuleshov and Shapiro~\cite{ks} for $d=3,4$ and $5$
can be obtained and extended, using
{\sf Maple} to carry out the routine manipulation
of series.

%

As an example, we give an explicit expression for
$H^{\genus{g}}_{(\al_1,\al_2),(\be_1,\be_2)}$, in the case
that $\al_1,\al_2,\be_1,\be_2$ are distinct. \lremind{twotwo}

\begin{corollary} \label{twotwo}
Let $\alpha=(\alpha_1,\alpha_2)\vdash d$ and
$\beta=(\beta_1,\beta_2)\vdash d$
where $\alpha_1<\alpha_2$, $\beta_1<\beta_2$, $\alpha_1<\beta_1$ and
$\alpha_1,\alpha_2,\beta_1,\beta_2$ are distinct.
Then
$$
H^{\genus{g}}_{\alpha,\beta} =
\frac{2}{\alpha_1\alpha_2\beta_1\beta_2}
\sum_{i=1}^{\alpha_1}\left(
\left(\binom {d+1}{ 2}-di\right) ^{2g+2}-
\left(\binom {d+1}{2}-di-\alpha_2\beta_1\right) ^{2g+2}
\right) .
$$
\end{corollary}

\proof
Since $\al_1,\al_2,\be_1,\be_2$ are distinct, we have
$$H^{\genus{g}}_{\alpha,\beta}=\tH^{\, (g)}_{\alpha,\beta}
=  \frac{\left|{\mathcal{C}}_{\beta}\right| \cdot\left|{\mathcal{C}}_{\alpha}\right| }{d!^2}
\sum_{\lambda\vdash d}\eta(\lambda )^r\chi^{\lambda}_{\alpha}
\chi^{\lambda}_{\beta},$$
from (\ref{discfactn}). But here we have $r=2g+2$ is even, and
$\left|{\mathcal{C}}_{\be}\right|
\cdot\left|{\mathcal{C}}_{\alpha}\right| /d!^2$ $=1/\al_1\al_2\be_1\be_2$. Moreover,
  $\eta(\widetilde{\la})=-\eta(\la )$, and $\chi^{\widetilde{\lambda}}_{\alpha}
\chi^{\widetilde{\lambda}}_{\beta}=\chi^{\lambda}_{\alpha}
\chi^{\lambda}_{\beta}$. Finally, from~(\ref{spps}),
$\chi^{\lambda}_{\alpha}= 0$
exactly when $[p_{\al}]s_{\la}=0$ , so we have
$$H^{\genus{g}}_{\alpha,\beta}=\frac{2}{\alpha_1\alpha_2\beta_1\beta_2}
\sum_{\lambda\in{\mathcal P}_{\al,\be}}\eta(\lambda)^{2g+2}
\chi^{\lambda}_{\alpha}\chi^{\lambda}_{\beta}, $$
where ${\mathcal P}_{\al,\be}=\{\la\colon
[p_{\al}]s_{\la}\neq0,[p_{\be}]s_{\la}\neq0,
\eta(\la)>0\}$.
Now, we can give an explicit description of ${\mathcal P}_{\al,\be}$, using
the Murnaghan--Nakayama formula
for the irreducible characters of the symmetric group (see, {\em e.g.,}
~\cite{macd}). The details can be routinely
verified. First,  $\left|{\mathcal P}_{\al,\be}\right| =2\alpha_1$, so we let
  ${\mathcal
P}_{\al,\be}=\left\{\lambda^{(1)},\ldots,\lambda^{(2\al_1)}\right\}$,
where $\lambda^{(1)}\succ\ldots\succ\lambda^{(2\al_1)}$.
(Here $\succ$ denotes reverse lexicographic order on partitions,
so $(3)\succ(2,1)\succ(1^3)$.)
Then, for $i=1,\ldots,\al_1$, we have
$\lambda^{(i)}=\left(d-i+1,1^{i-1}\right)$
(which is independent of $\beta$), so
$$\eta \left( \lambda^{(i)} \right) ={\binom {d+1} 2}-di,\;\;\;\;\;\;\;\;
\chi^{\lambda^{(i)}}_{\alpha}\chi^{\lambda^{(i)}}_{\beta}=1,\;\;\;\;
  i=1,\ldots,\alpha_1.$$
Also, for $i=1,\ldots,\al_1$, we have $\lambda^{(\al_1+i)}=\left(
d+1-\be_1-i,\al_1+2-i,2^{i-1},1^{\be_1-\al_1-1}
\right)$, so
$$\eta \left( \lambda^{(\al_1+i)} \right) ={\binom {d+1}  2}-di-\alpha_2\beta_1,\;\;\;\;\;\;\;\;
\chi^{\lambda^{(\al_1+i)}}_{\alpha}\chi^{\lambda^{(\al_1+i)}}_{\beta}=-1,\;\;\;\;
  i=1,\ldots,\alpha_1.$$
The result follows immediately.
\qed

For example,
${\mathcal P}_{(3,8),(4,7)}
=\left((11),(10,1),(9,1^2),(7,4),(6,3,2),(5,2^3)\right)_\succ$ and
$$
H^{\genus{g}}_{(3,8),(4,7)} =\frac{2}{3\cdot8\cdot4\cdot7}
\left(
55^{2g+2}+44^{2g+2}+33^{2g+2}-23^{2g+2}-12^{2g+2}-1^{2g+2}
\right).
$$
Similar expressions may be obtained when $\al$ and $\be$ have three parts.
For example,
\begin{eqnarray*}
H^{\genus{g}}_{(1,2,6),(1,3,5)} = \frac{1}{180}
\left(
2^{2g+4} - 6^{2g+4} + 10^{2g+4} + 12^{2g+4} -18^{2g+4} - 20^{2g+4}
-28^{2g+4} +36^{2g+4}
\right).
\end{eqnarray*}
The sum is over ${\mathcal P}_{\al,\be}$, but contributions from
some partitions of this set
are exactly canceled as a consequence of ``identities''
between parts of $\al$ and parts of $\be$ (for example, $1+2=3$ and $6=1+5$,
where the left and right hand sides, respectively, refer to $\al$ and $\be$).
Furthermore, other terms are introduced as a consequence of the same
identities.

As an example with more parts, with $d=8$, we obtain
\begin{eqnarray*}
H^{\genus{g}}_{(2,2,4),(1,2,2,3)} = \frac{1}{48}
\left(
 3 \cdot 2^{2g+5} + \frac 9 2 \cdot 4^{2g+5} +3\cdot 6^{2g+5}
- 10^{2g+5}   - 14^{2g+5} -16^{2g+5}
+ \frac 1 2 \cdot28^{2g+5}
\right).
\end{eqnarray*}

\section{$m$-part double Hurwitz numbers ($m=l(\al)$ fixed)}\label{marb}

\lremind{marb}We next consider more generally the case where $\al$ has a fixed
number $m$ of parts, and $\be$ is arbitrary.  The behaviour
is qualitatively different from that of the $m=1$ case, which was considered
in Section~\ref{m1}, as might be
expected by the failure of polynomiality.
We prove a
topological recursion relation consistent with a description of double
Hurwitz numbers in terms of the moduli space of curves.  This
relation is obvious neither from the currently understood geometry of
double Hurwitz numbers nor from the combinatorial interpretation in
terms of the join-and-cut equation.   For expository reasons, we
will give three versions of this topological recursion:  a 
genus $0$ recursion (Thm.\ ~\ref{gensymm}), a ``cleaner'' version
of the genus $0$ recursion involving rational rather than transcendental
functions
(Thm.\ ~\ref{eqminfour}), and a version
in arbitrary genus (Thm.\ ~\ref{intgenusg}).   

The topological recursion will enable us to find closed-form expressions
for double Hurwitz numbers for small $g$ and $m$, and in principle for
larger $g$ and $m$.  (For a much simpler example of topological recursions
implying closed-form expressions for single Hurwitz numbers, see \cite{v}.)
We conjecture the form of a closed-form expression for $g=0$ and 
arbitrary $m$ (Conj.\ ~\ref{conjecture1}).  

The reader will notice that except for the cases $(g,m)=(0,1)$ and $(0,2)$
(when there is no Deligne-Mumford moduli stack $\cmbar_{g,m}$), 
the explicit expressions that we obtain for $\HH_m^{\genus{g}}$ are
all rational functions in the intrinsic variable $u$. Moreover, the
denominator has explicit linear factors. The topological recursions
that we obtain for $\HH_m^{\genus{g}}$ are integrals over
$u$, and the integrand is quadratic in lower order terms. We
conjecture that $\HH_m^{\genus{g}}$ is, except for the two initial
cases, always a rational function in $u$, with specified linear factors
in the denominator. To prove this by induction, we would need to obtain a
rational integrand by the induction hypothesis, and then 
prove (to avoid a logarithm in the integrated form) that the inverse
linear terms in the partial fraction expansion of the integrand
disappear.  We have been unable to prove this in general, since it
seems to require a stronger induction hypothesis.

\subsection{The symmetrized join-cut equation at genus $0$} 
We apply the symmetrization operator $\Theta_m$ (defined in~(\ref{elso}))
to the join-cut equation \eqref{joincut} to obtain
partial differential equations for the symmetrized series
$\HH^{\genus{0}}_m(x_1,\ldots ,x_m)$.  
As a preliminary, we begin with $\HH^{\genus{0}}_1$.
The more general results will be an extension of this idea.

\begin{lemma}\label{eqonep} \lremind{eqonep}
$$\HH^{\genus{0}}_{1,1}(x_1)=uQ_1.$$ 
\end{lemma}

(Recall that $\HH^{\genus{g}}_{j,i}=x_i \frac {\partial \HH^{\genus{g}}_j} 
{ \partial x_i}.$)
Although Lemma~\ref{eqonep} has already been proved in the previous
section, in~(\ref{ans0}), we give a second proof to illustrate
the methodology that will be used throughout this section.

\proof
By applying $\Theta_1$ to the join-cut equation \eqref{joincut} and
setting $y=0$, it follows immediately that
$\HH^{\genus{0}}_1$ satisfies the partial differential equation
$$
\left( 1 + u \frac { \partial} {\partial u} + 0-2 \right)\HH^{\genus{0}}_1
= \left( u\frac{\partial}{\partial u} -1\right) \HH^{\genus{0}}_1
= \frac{1}{2} \left(  \left({\HH^{\genus{0}}}_{1,1} \right) ^2 + 0 + 0 \right)
=\frac{1}{2} \left({\HH^{\genus{0}}}_{1,1} \right)^2
$$
with initial condition $\left[ u \right]\HH^{\genus{0}}_{1,1} =Q(x_1)$.
Apply $\frac{x_1}{u^2}\frac{\partial}{\partial x_1}$ to the above equation
and let $G=\frac{1}{u}\HH^{\genus{0}}_{1,1}$, to obtain
\begin{equation} \label{diffsymone}
\frac{\partial G}{\partial u}=Gx_1\frac{\partial G}
{\partial x_1}.
\end{equation}
In terms of $G$, the initial condition becomes $\left[ u^0 \right]G=Q(x_1 )$.
But, applying $\frac{\partial}{\partial u}$ to the functional
equation~\eqref{funceqw}, we obtain \lremind{lagpartu}
\begin{equation}\label{lagpartu}
\frac{\partial w}{\partial u}=w\mu (w)Q(w),
\end{equation}
and comparing with~(\ref{lagpartial}), we check that $G(x_1)=Q_1$ is the
unique solution to~(\ref{diffsymone}).
\qed

To state an equation for $\HH^{\genus{0}}_m$ for
$m\geq 2$, we need some additional notation.  For $\alpha =\{\alpha
_1,\ldots ,\alpha _j\}\subseteq\{ 1,\ldots ,m\}$, let
$x_{\alpha}=x_{\alpha_1},\ldots ,x_{\alpha_j}$.  Let
$\Omega_{m,i}$ be the set of unordered pairs
$\{\alpha ,\zeta\}$ such that $\alpha ,\zeta\subseteq\{ 1,\ldots ,m
\}$ with $\alpha\cup\zeta =\{ 1,\ldots ,m\}$ and $\alpha\cap\zeta =\{
i\}$.  Let $\overline { \{ l \} } = \{ 1, \dots, m \} \setminus \{ l \}$.

\begin{theorem}[Symmetrized join-cut equation in genus $0$]\label{eqmp}
For $m\geq 2$, $\HH^{\genus{0}}_m$ satisfies the equation
\begin{eqnarray*}
\left( u\frac{\partial}{\partial u}+m-2
-\sum_{i=1}^m uQ_i
x_i\frac{\partial}{\partial x_i}\right) \HH^{\genus{0}}_m = \sum_{i=1}^m
\sum_{ \subst{\{\alpha ,\zeta\}\in \Omega_{m,i}}
{ l(\al), l(\zeta) \geq 2}  }
\!\!\!\!
\!\!\!\!
\HH^{\genus{0}}_{\left|\alpha \right| ,i}(x_{\alpha})
\HH^{\genus{0}}_{\left|\zeta \right| ,i}(x_{\zeta})
+\sum_{1\leq k,l\leq m\atop{k\neq l}}
\!\!\!\!
\frac{x_l \HH^{\genus{0}}_{m-1,k}(x_{\overline{\{ l\}}})}
{x_k -x_l}.
\end{eqnarray*}
\end{theorem}
The two parts of the right hand side of the above equation correspond to the
first two parts of the right hand side of the join-cut equation (Lemma~\ref{jc});
the third part of join-cut does not arise in genus $0$.

\proof
By applying $\Theta_m$ to~(\ref{joincut}) for fixed $m\geq 2$ and
setting $y=0$, we find
that $\HH^{\genus{0}}_m(x_1 ,\ldots ,x_m )$ satisfies
\begin{eqnarray*}
\left( u\frac{\partial}{\partial u}+m-2\right) \HH^{\genus{0}}_m =\sum_{i=1}^m
\sum_{\{\alpha ,\zeta\}\in \Omega_{m,i}}
\HH^{\genus{0}}_{\left|\alpha \right| ,i}(x_{\alpha}) \HH^{\genus{0}}_{\left|\zeta \right| ,i}(x_{\zeta})
+\sum_{1\leq k,l\leq m\atop{k\neq l}}
\!\!\!\!
\frac{x_l
\HH^{\genus{0}}_{m-1,k}(x_{\overline{\{ l\}}})}
{x_k -x_l}.
\end{eqnarray*}
Moving the contribution of $\{ \al, \zeta \} \in  \Omega_{m,i}$ where
$l(\al)=1$ or $l(\zeta)=1$
on the right hand side of this equation to the left hand side, we obtain
\begin{eqnarray*}
\left( u\frac{\partial}{\partial u}+m-2
-\sum_{i=1}^m \HH^{\genus{0}}_{1,i}(x_i)
x_i\frac{\partial}{\partial x_i}\right) \HH^{\genus{0}}_m = \sum_{i=1}^m
\sum_{ \subst{\{\alpha ,\zeta\}\in \Omega_{m,i}}
{ l(\al), l(\zeta) \geq 2}  }
\!\!\!\!
\!\!\!\!
\HH^{\genus{0}}_{\left|\alpha \right| ,i}(x_{\alpha})
\HH^{\genus{0}}_{\left|\zeta \right| ,i}(x_{\zeta})
+\sum_{1\leq k,l\leq m\atop{k\neq l}}
\!\!\!\!
\frac{x_l
\HH^{\genus{0}}_{m-1,k}(x_{\overline{\{ l\}}})}
{x_k -x_l},
\end{eqnarray*}
and the result follows from Lemma~\ref{eqonep}.
\qed

A key observation is the following. The right hand side of the
equation in Theorem~\ref{eqmp} involves the series $\HH^{\genus{0}}_j$ for $j<m$
only, so if we can invert the partial differential operator that is
applied to $\HH^{\genus{0}}_m$ on the left hand side, then we have a recursive
solution for $\HH^{\genus{0}}_m$, $m\geq 2$.

\subsection{Transformation of variables and recursive solution
to symmetrized join-cut}\label{ss:Grss} \lremind{ss:Grss}We 
now find a solution to the partial differential equation
for $\HH^{\genus{0}}_m$ that is given in Theorem~\ref{eqmp}.
  The key is
to change variables in $\HH^{\genus{0}}_m$,
for $m\geq 1$, from $x_1,\ldots ,x_m$ to
       $w_1,\ldots ,w_m$, using~(\ref{funceqw}), to obtain
\begin{equation}\label{changevars}
\hh^{\genus{0}}_m(u,w_1,\ldots ,w_m):=\HH^{\genus{0}}_m(w_1e^{-uQ_1},
\ldots ,w_me^{-uQ_m}),\;\;\;\;\;\;\;\; m\geq 2.
\end{equation}
We denote this transformation by $\Gamma$, so
$$\Gamma \HH^{\genus{0}}_m(x_1,\ldots,x_m)=\hh^{\genus{0}}_m(u,w_1,\ldots,w_m).$$
We regard $\hh^{\genus{0}}_m$ as an element of the ring of
formal power series in $u,w_1,\ldots ,w_m$, with
coefficients that are polynomials in $q_1,q_2,\ldots$.
It is straightforward to invert this, and
recover $\HH^{\genus{0}}_m$ from $\hh^{\genus{0}}_m$ in~(\ref{changevars})
by Lagrange inversion,
as specified in Theorem~\ref{lagthm}.
For this ring, let $D_u$ be the first partial derivative, with
respect to $u$, for the purposes of which $w_1,\ldots ,w_m$ are regarded as
algebraically independent variables, with no dependence on $u$.
Henceforth, we
use $\hh^{\genus{0}}_m$ and $\HH^{\genus{0}}_m$ interchangeably.

The importance of  $\Gamma$ is shown in its action on the partial
differential operator that is applied to $\HH^{\genus{0}}_m$
on the left hand side of the symmetrized join-cut
equation (Thm.~\ref{eqmp}).
We show that, under $\Gamma$, the partial differential
operator is transformed into a
linear differential operator involving only $D_u.$

\begin{lemma}\label{genint} \lremind{genint}Let $k$ be an integer. Then
$$ u^{k-1}
\Gamma\left( u\frac{\partial}{\partial u}+k
-\sum_{i=1}^m uQ_i
x_i\frac{\partial}{\partial x_i}\right)  = D_u u^k\Gamma.  $$
\end{lemma}

In short, passing $\Gamma$ through the differential
operator simplifies it.
  From this point of view, the variable $u$ plays an
important role, as the only variable, and accounts for our terming
it the intrinsic variable of the system.

\proof
For functions of $u,w_1,\ldots ,w_m$, the chain rule gives
\begin{equation*}
u \frac{\partial}{\partial u} = uD_u+\sum_{i=1}^m
u\left( \frac{\partial w_i}{\partial u}\right)
\frac{\partial}{\partial w_i}
       = uD_u+\sum_{i=1}^m uw_iQ_i\mu _i\frac{\partial}{\partial w_i}
= uD_u+\sum_{i=1}^m u Q_ix_i\frac{\partial}{\partial x_i},
\end{equation*}
from~(\ref{lagpartu}) and the operator identity~(\ref{opident}).
Then
$
\Gamma\left( u\frac{\partial}{\partial u}+k
-\sum_{i=1}^m uQ_i
x_i\frac{\partial}{\partial x_i}\right)  = \left(uD_u +k \right)\Gamma
$
and the result follows. \qed

\subsection{A (univariate, rational, integral) topological recursion
for $\HH_m^{\genus{0}}$}

Lemma~\ref{genint} enables us to solve the partial differential equation
for $\HH^{\genus{0}}_m$, $m\geq 2$, recursively.  The following result gives an
integral expression for $\hh^{\genus{0}}_m$, in terms of
$\hh^{\genus{0}}_1,\ldots ,\hh^{\genus{0}}_{m-1}$. We use the notation
$$\hh^{\genus{0}}_{j,i}=w_i\frac{\partial \hh^{\genus{0}}_j}{\partial w_i}.$$  
(A ``cleaner'' version, not involving
exponentials, will be given later, Theorem~\ref{eqminfour}.)

\begin{theorem}[Genus $0$ topological recursion, transcendental 
form]\label{gensymm}
\lremind{gensymm}For $m\geq 2$,
\begin{eqnarray*}
\hh^{\genus{0}}_m
= u^{2-m}\int_0^u\left(\sum_{i=1}^m
\sum_{ \subst{\{\alpha ,\zeta\}\in \Omega_{m,i}}
{ l(\al), l(\zeta) \geq 2}  }
\!\!\!\!
\!\!\!\!
\mu_i^2 \hh^{\genus{0}}_{\left|\alpha \right| ,i}(u,w_{\alpha})
\hh^{\genus{0}}_{\left|\zeta \right| ,i}(u,w_{\zeta})
+\sum_{1\leq k,l\leq m\atop{k\neq l}}
\!\!\!\!
\frac{w_le^{-uQ_l}\mu_k
\hh^{\genus{0}}_{m-1,k}(u,w_{\overline{\{ l\}}})}
{w_ke^{-uQ_k} -w_le^{-uQ_l}}\right)u^{m-3}du,
\end{eqnarray*}
where the integrand is considered as a power series in $u,w_1,\ldots ,w_m$, and
the integration is carried out with $w_1,\ldots ,w_m$ regarded as constants.
\end{theorem}

\proof
The result follows by applying $\Gamma$ to Theorem~\ref{eqmp} with the
aid of Lemma~\ref{genint}.
\qed

\subsection{Explicit expressions for $\HH_m^{\genus{0}}$ for $m\leq 5$, and
a conjectured form in general}
We now apply Theorem~\ref{gensymm} for successive values
of $m\geq 2$, to obtain explicit expressions for the symmetrized
series $\HH_m^{\genus{0}}$. We begin with $m=2$ and $m=3$, and
include the details in a single result, because the resulting
expressions can be treated uniformly, also incorporating $m=1$.
This requires some notation.
For $m\geq 1$, let $V_m=\prod_{1\leq i< j\leq m}(w_i-w_j)$, the value of the
Vandermonde determinant $\det\left( w_j^{m-i}\right) _{m\times m}$,
and let $\bA_m$ be the $m\times m$ matrix with $(1,j)$-entry equal
to $\mu _j-1$, for $j=1,\ldots ,m$, and $(i,j)$-entry equal
to $w_j^{m-i+1}$, for $i=2,\ldots ,m,\;\; j=1,\ldots ,m$.
Let $\Delta_{m,j}$ be
the partial differential operator defined by
\begin{equation}\label{deltadef}
\Delta_{m,j}=\sum_{i=1}^m w_i^j\mu _i\frac{\partial}{\partial w_i}
=\sum_{i=1}^m w_i^{j-1}x_i\frac{\partial}{\partial x_i},
\end{equation}
for $m,j\geq 1$, where the second equality follows from the
operator identity~(\ref{opident}).

\begin{corollary}\label{twoparts}
\lremind{twoparts}For $m=1,2,3$,
$$\Delta_{m,1}^{3-m}\hh^\genus{0}_m=\frac{\det \bA_m}{V_m}.$$
\end{corollary}

\proof
For $m=1$, the result follows by differentiating the result of
Lemma~\ref{eqonep} and applying~(\ref{lagpartial}).

For the case $m=2$, Theorem~\ref{gensymm} gives

\begin{eqnarray*}
\hh^{\genus{0}}_2
&=& \int_{0}^u\left(\frac{uQ_1w_2e^{-uQ_2}-uQ_2w_1e^{-uQ_1}}
{w_1e^{-uQ_1}-w_2e^{-uQ_2}}-Q_1-Q_2\right) \frac{du}{u} \\
&=& \int_{0}^u\left(\frac{Q_1w_1e^{-uQ_1}-Q_2w_2e^{-uQ_2}}
{w_1e^{-uQ_1}-w_2e^{-uQ_2}}-Q_1-Q_2\right) du \\
&=&-\log\left( \frac{w_1e^{-uQ_1}-w_2e^{-uQ_2}}{w_1-w_2}\right)
-uQ_1-uQ_2 \\
&=& \log\left( \frac{w_1-w_2}{x_1-x_2}\right) -\left( uQ_1+uQ_2\right),
\end{eqnarray*}
which is well-formed as a formal power series in $u,w_1,w_2$, since
the argument for the logarithm has constant term equal to $1$.
Then from~(\ref{lagpartial}) we have
\begin{equation}\label{e12}
\HH^{\genus{0}}_{2,1}=E_{1,\, 2}-\frac{x_2}{x_1-x_2},\;\;\;\,
\HH^{\genus{0}}_{2,2}=E_{2,\, 1}-\frac{x_1}{x_2-x_1},
\quad\mbox{where}\quad
E_{i,\, j}=\frac{w_j\mu _i}{w_i -w_j},
\end{equation}
for $i\neq j$, and, adding these, we obtain
\begin{eqnarray*}
\left( x_1\frac{\partial}{\partial x_1}+x_2\frac{\partial}{\partial x_2}
\right) \HH^{\genus{0}}_2 &=& E_{1,\, 2}+E_{2,\, 1}+1
= \frac{w_2\mu _1}{w_1-w_2}
+\frac{w_1\mu _2}{w_2-w_1}+1 \\
&=&\frac{w_2}{w_1-w_2}(\mu _1-1)+
\frac{w_1}{w_2-w_1}(\mu _2-1)
\end{eqnarray*}
so the result follows for $m=2$.

For $m=3$,
let $x_{i,j}=x_j/(x_i-x_j)$, for $i\neq j$, and let $\sum_{i,j,k}$
       denote summation over all distinct $i,j,k$ with $1\leq i,j,k\leq 3$.
For the case $m=3$,
Theorem~\ref{gensymm} and~(\ref{lagpartial}),~(\ref{e12}) give
\begin{eqnarray*}
\hh^{\genus{0}}_3
&=&\frac{1}{u}\int_0^u \sum_{i,j,k}\left(\frac12
(E_{i,j}-x_{i,j})(E_{i,k}-x_{i,k})
       +x_{i,j}(E_{i,k}-x_{i,k})\right) du \\
&=&\frac{1}{u}\int_0^u\left( -1+ E_{1, 2}E_{1, 3}+E_{2, 1}E_{2, 3}+
E_{3, 1}E_{3, 2}
\right) du \\
&=&\frac{1}{u}\int_0^u\left( -1+\sum_{i=1}^3 \left(\prod_{1\leq j\leq
3\atop{j\neq i}}
\frac{w_j}{w_i-w_j}\right)\mu_i^2\right) du
\\
&=&\frac{1}{u}\left. \left( -u+\sum_{i=1}^3
\left(\prod_{1\leq j\leq 3\atop{j\neq i}}
\frac{w_j}{w_i-w_j}\right)\frac{\mu_i}
{w_iQ_i^{\prime}}
\right) \right| ^u_0\\
&=&-1+\sum_{i=1}^3
\left(\prod_{1\leq j\leq 3\atop{j\neq i}}
\frac{w_j}{w_i-w_j}\right) \mu _i 
= \sum_{i=1}^3
\left(\prod_{1\leq j\leq 3\atop{j\neq i}}
\frac{w_j}{w_i-w_j}\right) (\mu _i-1),
\end{eqnarray*}
and the result follows for $m=3$.
\qed


We now apply the case $m=2$ of the previous result to transform
the integrand of Theorem~\ref{gensymm} to a simpler form.
\begin{theorem}[Genus $0$ topological recursion, rational form]\label{eqminfour}
\lremind{eqminfour}For $m\geq 4$,
\begin{equation*}
\hh^{\genus{0}}_m
=u^{2-m}\int_0^u\left(\sum_{i=1}^m
\sum_{ \subst{\{\alpha ,\zeta\}\in \Omega_{m,i}}
{ l(\al), l(\zeta) \geq 3}  }
\mu_i^2 \hh^{\genus{0}}_{\left|\alpha \right| ,i}(u,w_{\alpha})
\hh^{\genus{0}}_{\left|\zeta \right| ,i}(u,w_{\zeta})
+\sum_{1\leq k,l\leq m\atop{k\neq l}}\frac{\mu_k^2w_l}
{w_k-w_l}
\hh^{\genus{0}}_{m-1,k}(u,w_{\overline{\{ l\}}})\right)u^{m-3}du.
\end{equation*}
\end{theorem}

\proof
The result follows immediately from Theorem~\ref{gensymm}
and~(\ref{e12}), since we are able to cancel the terms with
denominator $x_k-x_l=w_ke^{-uQ_k}-w_le^{-uQ_l}$.
\qed

We now apply Theorem~\ref{eqminfour} in the cases
of $m=4$ and $m=5$ parts.
For $m\geq 3$, let $\bB_m^{(n;k)}$ be the $m\times m$ matrix
with $(1,j)$-entry equal to $w_j^n\mu_j$, $(2,j)$-entry
equal to $w_j\mu_j$, $(3,j)$-entry equal to $w_j^k$,
for $j=1,\ldots m$,
and $(i,j)$-entry equal to $w_j^{m-i+1}$,
for $i=4,\ldots ,m,\;\; j=1,\ldots ,m$.\lremind{fourparts}

\begin{corollary}\label{fourparts}
\begin{equation*}
\hh^{\genus{0}}_4=\Delta_{4,1}\left(\frac{\det \bA_4}{V_4} \right)
-\frac{\det \bB_4^{(2;2)}}{V_4}.
\end{equation*}
\end{corollary}

The proof is similar in approach to that of Corollary~\ref{twoparts}
and therefore omitted.

Note that the right hand side of Corollary~\ref{twoparts} with $m=4$
appears as the
``first approximation'' to the series $\hh^{\genus{0}}_4$ in Corollary~\ref{fourparts}.

For $m=5$, we have found the expressions in Theorem~\ref{eqminfour} to be
intractable by hand, but have used {\sf Maple} to carry out
the integration, and obtained the following theorem.
Let
\begin{equation*}
\bC_5=
\begin{pmatrix}
w_1^4\mu_1 & w_2^4\mu_2 & w_3^4\mu_3 & w_4^4\mu_4 & w_5^4\mu_5 \cr
w_1^2\mu_1 & w_2^2\mu_2 & w_3^2\mu_3 & w_4^2\mu_4 & w_5^2\mu_5 \cr
w_1\mu_1 & w_2\mu_2 & w_3\mu_3 & w_4\mu_4 & w_5\mu_5 \cr
w_1^2 & w_2^2 & w_3^2 & w_4^2 & w_5^2 \cr
w_1 & w_2 & w_3 & w_4 & w_5 \cr
\end{pmatrix}
,\quad
\bD_5^{(n;k)}=
\begin{pmatrix}
w_1^n\mu_1^2 & w_2^n\mu_2^2 & w_3^n\mu_3^2 & w_4^n\mu_4^2 & w_5^n\mu_5^2 \cr
w_1^k\mu_1 & w_2^k\mu_2 & w_3^k\mu_3 & w_4^k\mu_4 & w_5^k\mu_5 \cr
w_1^3 & w_2^3 & w_3^3 & w_4^3 &w_5^3 \cr
w_1^2 & w_2^2 & w_3^2 & w_4^2 & w_5^2 \cr
w_1 & w_2 & w_3 & w_4 & w_5 \cr
\end{pmatrix}.
\end{equation*}

\begin{corollary}\label{fiveparts}
\begin{eqnarray*}
\hh^{\genus{0}}_5 &=& \Delta _{5,1}^2\left(\frac{\det \bA_5}{V_5} \right)
-\Delta _{5,1}\left( 2\frac{\det \bB_5^{(3;3)}}{V_5}-\frac{\det
\bB_5^{(2;4)}}{V_5}\right)
+\Delta_{5,2}\left(\frac{\det \bB_5^{(2;3)}}{V_5}\right)  \\
&\mbox{}&  +\frac{\det \bC_5}{V_5}
+\frac{\det \bD_5^{(3;1)}}{V_5}
-\frac{\det \bD_5^{(2;2)}}{V_5}.
\end{eqnarray*}
\end{corollary}

Again, the right hand side of Corollary~\ref{twoparts} with $m=5$
appears as the ``first approximation'' to $\hh^{\genus{0}}_5$ in
Corollary~\ref{fiveparts}.

The results that we have for $m=1, \dots, 5$ have not yet suggested a
pattern that can be conjecturally generalised.  This is because we
have been unable to find a sufficiently uniform presentation for them,
although the presentation as a sum of bialternants of very elementary matrices
is appealing.

Still, the forms that we have obtained for $\hh^{\genus{0}}_m$ when $m\leq 5$ suggest a
general conjecture, stated below. We refer to this as a rational
form in $u$, because each $\mu_i$ is an inverse linear function
of $u$. Note that, for $\hh^{\genus{0}}_m$ to continue to be rational as $m$
   increases, the partial fraction expansion of the recursively
formed integrand in Theorem~\ref{eqminfour} must continue to have
vanishing coefficients for the terms that are linear in $\mu_i$,
$i=1,\ldots ,m$. 

\begin{conjecture}\label{conjecture1}
\lremind{conjecture1}For $m\geq 3$,
$\hh^{\genus{0}}_m$ is a sum of terms of the following type:
$$\Delta_{m,i_1}\ldots\Delta_{m,i_k}P_{m,i_1,\ldots ,i_k},$$
where $0\leq k\leq m-3$, $i_1+ \ldots +i_k\leq m-3$, and
$P_{m,i_1,\ldots ,i_k}$
       is a homogeneous symmetric polynomial in $\mu_1,\ldots ,\mu_m$ of
degree $k+1$, with coefficients that are rational functions
in $w_1,\ldots ,w_m$ with degree of numerator minus degree of
denominator equal to $i_1 + \cdots +i_k -k$. Moreover,
       $P_{m,i_1,\ldots ,i_k}$ is a symmetric function of $w_1,\ldots ,w_m$,
where $\mu_i$ is considered as $\mu(w_i)$.
\end{conjecture}

Note that this form specializes to the expressions above 
for $\hh^{\genus{0}}_3, \hh^{\genus{0}}_4, \hh^{\genus{0}}_5$,
so the conjecture is true for the cases $3\leq m\leq 5$.

This conjecture should be seen as the genus $0$ double Hurwitz analogue
of the polynomiality conjecture \cite[Conj.\ 1.2]{gj2} (proved in
\cite[Thm.\ 3.2]{gjvk}).  As with the earlier conjecture, the form of
Conjecture~\ref{conjecture1} suggests some geometry.  For example,
present in the polynomial conjecture was the dimension of the moduli
space of $n$-pointed genus $g$ curves; the $n$ points corresponded to
the preimages of $\infty$.  In this case, the analogue is $m-3$, the
dimension of the moduli space of $m$-pointed genus $0$ curves; again,
the $m$ points should correspond to the preimages of $\infty$ ({\em i.e.}
the parts of $\beta$).  However, we have been unable to make precise
the link to geometry.

\subsection{Application:  explicit formulae}\label{m45}
As an application of the explicit formulae for $\hh^{\genus{0}}_m$ for
small $m$, we now extract the appropriate coefficient to
give explicit formulae for the corresponding double Hurwitz numbers.
We use some standard results for symmetric functions (see, for
example~\cite{macd}), particularly the
determinantal identity
\begin{equation}\label{jactrud}
\frac{\det \left(
w_{j}^{\theta _{i}+m-i}\right) _{m\times m} }{V_{m}}
=\det \left( h_{\theta _{i}-i+j}(\mathbf{w})\right) _{m\times m} ,
\end{equation}
for non-negative integers $\theta _{1},\ldots ,\theta _{m}$,
where $h_{k}(\mathbf{w})$
is the complete symmetric function of total degree $k$, with generating
series $\sum_{k\geq 0} h_k(\mathbf{w})t^k=\prod_{j=1}^m(1-w_jt^j)^{-1}$.
If $\theta =(\theta_1,\ldots ,\theta_m )$ is a partition (where $\theta_1\geq
\ldots\geq\theta_m$), then both sides of~(\ref{jactrud}) give
expressions for the Schur symmetric function $s_{\theta}(\mathbf{w})$.
In the case that $\theta$ is not a partition, we shall still denote
either side of~(\ref{jactrud}) by $s_{\theta}(\mathbf{w})$.

Using multilinearity on the first row of $\det \bA_m$,
we have
\begin{equation}\label{Amvanx}
\frac{\det \bA_{m}}{V_{m}}= w_{1}\ldots
w_{m} \sum_{r\ge m}a_r s_{\left( r-m\right) }\left( \mathbf{w}\right)
=w_{1}\ldots w_{m} \sum_{r\ge m}a_rh_{ r-m }(\mathbf{w}),
\end{equation}
from~(\ref{jactrud}), where 
\begin{eqnarray}\label{ear}
\mu(w)=\sum_{i\ge1}a_iw^i,
\end{eqnarray}
and $\mu(w)$ is defined in~(\ref{lagpartial}).
We write $\al\cup\be$ for the partition with parts
$\al_1,\ldots ,\al_m, \be_1,\ldots ,\be_n$, suitably reordered.

\begin{proposition}\label{detAmexpn} \lremind{detAmexpn}
For $d\geq m\geq 1$ and $\al ,\be\vdash d$, with
  $\al=(\al_1,\ldots ,\al_m)$,
$$\left[ x_1^{\al_1}\ldots x_m^{\al_m}u^{l(\be )}q_{\be} \right]\frac{\det \bA_m}{V_m} =\sum
\frac{l(\rho )!\prod_{j\geq 1}\rho_j}{\left| \Aut \rho  \right| }
\prod_{j=1}^m\frac{(\al_j-\left|\gamma_j\right| )\al_j^{l(\gamma_j)-1}}
{\left|\Aut
\gamma_j \right| },
$$ 
where the summation  is over partitions $\rho,\gamma_1,\ldots ,\gamma_m$,
with $\rho\cup\gamma_1\cup\ldots\cup\gamma_m=\be$, and
$\left|\gamma_j\right| <\al_j$, $j=1, \dots ,m$.
(Note that $\gamma_1,\dots ,\gamma_m$ can be empty, but $\rho$ cannot.)
\end{proposition}

\proof From~\eqref{Amvanx}, we have
$$\frac{\det \bA_m}{V_m} = \sum_{r\geq m}a_r \sum w_1^{i_1}\cdots w_m^{i_m},$$
where $m\geq 1$, and  the second summation is over $i_1, \dots ,i_m\geq 1$,
with $i_1+\cdots +i_m=r$, and $a_r$ is defined above~\eqref{ear}.
But, applying Theorem~\ref{lagthm}\eqref{lagthm1} to~(\ref{funceqw}), we obtain
\begin{eqnarray*}
\left[ x^t \right]w^i&=&\frac{1}{t} \left[ \la^{t-1} \right] i\la^{i-1}e^{utQ(\la )}=
\frac{i}{t} \left[ \la^{t-i} \right]\sum_{n\geq 0}\frac{1}{n!}
\left( ut\sum_{j\geq 1}q_j\la ^j\right) ^n\\
&=& \frac{i}{t}\sum_{\gamma\vdash t-i} \frac{t^{l(\gamma )}}{ \left|\Aut \gamma  \right| }
u^{l(\gamma )} q_{\gamma},\;\;\;\; i,t\geq 1,
\end{eqnarray*}
and from~(\ref{lagpartial}),
we have
\begin{eqnarray*}
a_r&=& \left[ w^r \right]\left(\mu(w)-1\right) =\left[ w^r \right] \sum_{i\geq 1}\left(
u\sum_{j\geq 1}jq_jw^j\right) ^i\\
&=&\sum_{\rho\vdash r}\frac{l(\rho )!\prod_{j\geq 1}\rho_j}{ \left|\Aut \rho 
\right| }
u^{l(\rho )}q_{\rho},\;\;\;\; r\geq 1.
\end{eqnarray*}
Combining these results, we obtain
$$\left[ x_1^{\al_1}\ldots x_m^{\al_m} \right]\frac{\det \bA_m}{V_m} =\sum
\frac{l(\rho )!\prod_{j\geq 1}\rho _j}{ \left| \Aut \rho \right|  }
u^{l(\rho )}q_{\rho}
\prod_{j=1}^m \left(\frac{i_j}{\al _j}\sum_{\gamma _j\vdash\al _j-i_j}
\frac{{\al _j}^{l(\gamma _j)}}{\left| \Aut \gamma _j\right| }
u^{l(\gamma _j)} q_{\gamma _j} \right) ,$$
where the summation is over $\rho$ and $i_1,\ldots ,i_m\geq 1$,
with $\rho \vdash i_1+\cdots +i_m$.
The result follows immediately.
\qed

This result allows us to immediately give formulae for genus $0$ double
Hurwitz numbers when one of the partitions has two or three parts.

\begin{corollary}\label{coef23pts}
Suppose $\al ,\be\vdash d$, with $\al =(\al_1,\ldots ,\al_m)$.

1) If $m=2$, then
$$H^{\genus{0}}_{(\al _1,\al _2),\be}=\frac{ \Autbe r!}{d}
\sum\frac{l(\rho )!\prod_{j\geq 1}\rho_j}{\left| \Aut \rho \right|  
\left|\Aut \gamma _1\right| 
\left| \Aut \gamma_2\right| }
(\al_1-\left|\gamma_1\right| )(\al_2-\left|\gamma_2\right| )\al_1^{l(\gamma_1)-1}\al_2^{l(\gamma_2)-1},$$
where the summation  is over partitions $\rho,\gamma_1,\gamma_2$,
with $\rho\cup\gamma_1\cup\gamma_2=\be$, and $\left|\gamma_j\right| <\al_j$, $j=1,2$.

2) If $m=3$, then
$$H^{\genus{0}}_{(\al _1,\al _2,\al _3),\be}= \Autbe r! 
\sum\frac{l(\rho )!\prod_{j\geq 1}\rho_j}{ \left| \Aut \rho  \right| }
\prod_{j=1}^3\frac{(\al_j-\left|\gamma_j\right| )\al_j^{l(\gamma_j)-1}}{\left| \Aut
\gamma _j \right| },$$
where the summation  is over partitions $\rho,\gamma_1,\gamma_2,\gamma_3$,
with $\rho\cup\gamma_1\cup\gamma_2\cup\gamma_3=\be$, and
$\left|\gamma_j\right| <\al_j$, $j=1,2,3$.
\end{corollary}

\proof From~(\ref{symtransf}), we obtain
$$H^{\genus{0}}_{(\al _1,\ldots ,\al _m),\be}=
\Autbe r! 
\left[ x_1^{\al_1}\ldots x_m^{\al _m}u^{l(\be )}q_{\be} \right] \HH^{\genus{0}}_m(x_1,\ldots ,x_m).$$
Both parts of the result then follow from Proposition~\ref{detAmexpn} and
Corollary~\ref{twoparts}, using~(\ref{deltadef}) to give the
factor of~$d$ in the case $m=2$.
\qed

In a similar, but more complicated way, it is possible to obtain explicit
formulae for $H^{\genus{0}}_{(\al _1,\ldots ,\al _m),\be}$ in the
cases $m=4,5$, using multilinearity to expand the determinants
that arise in Corollaries~\ref{fourparts} and~\ref{fiveparts}. 

\subsection{Positive genus:  A topological recursion for $\HH_m^{\genus{g}}$ and
  explicit formulae} In the following result, we apply the
symmetrization operator $\Theta_m$ to the join-cut equation, to obtain
a partial differential equation for $\HH_m^{\genus{g}}$, for genus
$g\geq 1$. As in the case of genus $0$, the change of variables
transforms the partial differential operator applied to
$\HH_m^{\genus{g}}$ into the linear differential operator in the
intrinsic variable $u$. Consequently, we are able to express the
transformed series
$$\hh_m^{\genus{g}}(u,w_1,\ldots ,w_m)=\Gamma \HH_m^{\genus{g}}(x_1,\ldots,x_m)$$ as
an integral in $u$. \lremind{intgenusg}

\begin{theorem}[Topological recursion in positive genus]\label{intgenusg}
\quad

\noindent
1) For $g\geq 1$,

\begin{eqnarray*}
\hh_1^{\genus{g}}=\frac{u^{1-2g}}{2}\int_0^u\left(
\sum_{j=1}^{g-1}\hh^{\genus{j}}_{1,1}(u,w_1)\hh^{\genus{g-j}}_{1,1}(u,w_1)
+\left. w_{2}\frac{\partial}{\partial w_2}
\hh^{\genus{g-1}}_{2,1}(u,w_1,w_2)\right| _{w_2=w_1}\right)
\mu_1^2 u^{2g-2}du.
\end{eqnarray*}

\noindent
2) For $m\geq 2$ and $g\geq 1$,
\begin{eqnarray*}
\hh_m^{\genus{g}}=u^{2-m-2g}\int_0^u \left(
\sum_{i=1}^m
\sum_{ \subst{\{\alpha ,\zeta\}\in \Omega_{m,i}}
{ l(\al), l(\zeta) \geq 3}  }
\mu_i^2 \left(
\hh^{\genus{0}}_{\left|\alpha \right| ,i}(u,w_{\alpha})
\hh^{\genus{g}}_{\left|\zeta \right| ,i}(u,w_{\zeta})
+ \hh^{\genus{g}}_{\left|\alpha \right| ,i}(u,w_{\alpha})
\hh^{\genus{0}}_{\left|\zeta \right| ,i}(u,w_{\zeta})
\right)  \right.
\end{eqnarray*}
$$+\sum_{j=1}^{g-1}\sum_{i=1}^m
\sum_{\{\alpha ,\zeta\}\in \Omega_{m,i}}\mu_i^2
\hh^{\genus{j}}_{\left|\alpha \right| ,i}(u,w_{\alpha})
\hh^{\genus{g-j}}_{\left|\zeta \right| ,i}(u,w_{\zeta})
+\sum_{1\leq k,l\leq m\atop{k\neq l}}
\!\!\!\!
\frac{\mu_k^2w_l}{w_k-w_l}\hh^{\genus{g}}_{m-1,k}(u,w_{\overline{\{ l\}}})$$
$$+\left. \frac{1}{2}\sum_{i=1}^m\mu_i^2 \left.\left(
w_{m+1}\frac{\partial}{\partial w_{m+1}}
\hh^{\genus{g-1}}_{m+1,i}(u,w_1,\ldots ,w_{m+1}) \right) \right| _{w_{m+1}=w_i}
\right) u^{m+2g-3}du .$$
In both parts of this result, the integration is carried out with
$w_1,\ldots ,w_m$ regarded as constants.
\end{theorem}

We call this  a topological recursion because it expresses
$\hh^{\genus{g}}_m$ in terms of $\hh^{\genus{g'}}_{m'}$, where $g' \leq g$
and $m' \leq m+1$, and either $g' <g$ or $m'<m$.

\noindent
{\em Remarks.}
\newline \noindent
{\em 1.}
Note that the case $m=1$ is different.
\newline \noindent
{\em 2.}
The exponents of $u$ have geometric meaning; this
is no coincidence.
\newline \noindent
{\em 3.}
 This result specializes to the rough form of the genus
0 topological recursion (Theorem~\ref{gensymm}), by taking
$g=0$ and $\hh^{\genus{-1}}=0$, after minor manipulation.

\proof
By applying $\Theta_1$ and $\left[ y^g \right]$ to~(\ref{joincut}) for fixed
$g\geq 1$ we find that $\HH^{\genus{g}}_1(x_1)$ satisfies
\begin{eqnarray*}
\left( u\frac{\partial}{\partial u}+2g-1\right) \HH^{\genus{g}}_1(x_1) =
\frac{1}{2}\sum_{j=0}^{\genus{g}}
\HH^{\genus{j}}_{1,1}(x_1)
\HH^{\genus{g-j}}_{1,1}(x_1)
+\frac{1}{2}\left. x_{2}\frac{\partial}{\partial x_2}
\HH^{\genus{g-1}}_{2,1}(x_1,x_2)\right| _{x_2=x_1}.
\end{eqnarray*}
Now move the terms $j=0$ and $j=g$ in the summation on the right hand
side of this equation
to the left hand side, and change variables by applying the
operator identity~(\ref{opident}), and part 1 of the result
follows from Lemma~\ref{genint}.

By applying $\Theta_m$ and $\left[ y^g \right]$ to~(\ref{joincut}) for fixed $m\geq 2$ and
$g\geq 1$ we find that $\HH^{\genus{g}}_m(x_1 ,\ldots ,x_m )$ satisfies
\begin{eqnarray*}
\left( u\frac{\partial}{\partial u}+m+2g-2\right) \HH^{\genus{g}}_m &=&
\sum_{j=0}^{\genus{g}}\sum_{i=1}^m
\sum_{\{\alpha ,\zeta\}\in \Omega_{m,i}}
\HH^{\genus{j}}_{\left|\alpha \right| ,i}(x_{\alpha})
\HH^{\genus{g-j}}_{\left|\zeta \right| ,i}(x_{\zeta})
+\sum_{1\leq k,l\leq m\atop{k\neq l}}
\!\!\!\!
\frac{x_l
\HH^{\genus{g}}_{m-1,k}(x_{\overline{\{ l\}}})}
{x_k -x_l}\\
& & +\frac{1}{2}\sum_{i=1}^m\left.\left( x_{m+1}\frac{\partial}{\partial x_{m+1}}
\HH^{\genus{g-1}}_{m+1,i}(x_1,\ldots ,x_{m+1})\right)\right| _{x_{m+1}=x_i}.
\end{eqnarray*}

Now move the contribution of $\{ \al, \zeta \} \in \Omega_{m,i}$ when $j=0$ or $j=g$,
and $l(\al)=1$ or $l(\zeta)=1$, in the first summation on
the right hand side of this equation to the left hand side, and
apply~(\ref{e12}) to  cancel the terms with denominator $x_k-x_l$ on
the right hand side. Then apply $\Gamma$, using the operator
identity~(\ref{opident}), and part 2 of the result follows
from Lemma~\ref{genint}. \qed


The topological recursion may be used to give explicit formulae for
$\HH_m^{\genus{g}}$.  
The cases $g=1$ and $m=1,2$ are given below.  We omit the derivation (which is
similar in spirit to that for the genus $0$ formulae), and simply report
the result.

\begin{corollary} \label{ttp} \lremind{ttp}
\begin{eqnarray*}
\hh^{\genus{1}}_1 &=&
u \mu_1w_1\frac{\partial}{\partial w_1}Q_1 +
u^2\mu_1^3   \left(\left(   w_1\frac{\partial}{\partial w_1}\right)^2
Q_1\right)^2 +
u\mu_1^2 \left(   w_1\frac{\partial}{\partial w_1}\right)^3 Q_1. \\
\hh_2^{\genus{1}}&=& \frac{1}{24}
\left( x_1\frac{\partial}{\partial x_1}\right)^2
\left(
\frac{w_2}{w_1-w_2} w_1\frac{\partial}{\partial w_1} \log \mu_1
\right)
+ \frac{1}{24}
\left( x_2\frac{\partial}{\partial x_2}\right)^2
\left(
\frac{w_1}{w_2-w_1} w_2\frac{\partial}{\partial w_2} \log \mu_2
\right) \\
&\mbox{}& -\frac{1}{24} \Delta_{2,1}
\left(\frac{w_2}{w_1-w_2}\mu_1 + \frac{w_1}{w_2-w_1}\mu_2 \right)
   +
\frac{1}{48}  \Delta_{2,1}^2
\left(
w_1\frac{\partial}{\partial w_1}
\frac{w_2}{w_1-w_2}
+
w_2\frac{\partial}{\partial w_2}
\frac{w_1}{w_2-w_1}
\right).
\end{eqnarray*}
\end{corollary}

A similar equation for $\hh^\genus{1}_1$
has also been derived from the Genus Expansion Ansatz,
in~(\ref{ans1}). Of course, these expression agree, by carrying
out the differentiations in Corollary~\ref{ttp}.



\end{document}

%% file: cornereg4.tex
\setlength{\unitlength}{0.00083333in}
\begingroup\makeatletter\ifx\SetFigFont\undefined%
\gdef\SetFigFont#1#2#3#4#5{%
  \reset@font\fontsize{#1}{#2pt}%
  \fontfamily{#3}\fontseries{#4}\fontshape{#5}%
  \selectfont}%
\fi\endgroup%
{\renewcommand{\dashlinestretch}{30}
\begin{picture}(3624,2439)(0,-10)
\put(1699.500,1099.500){\arc{530.330}{4.5705}{6.4251}}
\path(12,462)(762,762)(1512,462)
	(2262,612)(2862,1062)(3612,462)
\path(12,1962)(612,1512)(762,762)
	(1662,1062)(2862,1062)(2712,1812)
	(1662,1962)(1662,1062)(612,1512)
	(1662,1962)(1512,2412)
\path(612,1512)(912,2412)
\path(2712,1812)(3462,2412)
\path(2862,1062)(3612,1512)
\path(2262,612)(3162,12)
\path(2262,612)(2112,12)
\path(1512,462)(1062,12)
\path(2187,1587)(1812,1212)
\blacken\path(1875.640,1318.066)(1812.000,1212.000)(1918.066,1275.640)(1871.397,1271.397)(1875.640,1318.066)
\put(2262,1662){\makebox(0,0)[lb]{\smash{{{\SetFigFont{8}{9.6}{\rmdefault}{\mddefault}{\updefault}corner}}}}}
\end{picture}
}

%% file: corneredge2.tex
\setlength{\unitlength}{0.00083333in}
\begingroup\makeatletter\ifx\SetFigFont\undefined%
\gdef\SetFigFont#1#2#3#4#5{%
  \reset@font\fontsize{#1}{#2pt}%
  \fontfamily{#3}\fontseries{#4}\fontshape{#5}%
  \selectfont}%
\fi\endgroup%
{\renewcommand{\dashlinestretch}{30}
\begin{picture}(3624,2439)(0,-10)
\put(1672,2092){\blacken\ellipse{50}{50}}
\put(1672,2092){\ellipse{50}{50}}
\put(1577,2062){\blacken\ellipse{50}{50}}
\put(1577,2062){\ellipse{50}{50}}
\put(1462,1932){\blacken\ellipse{50}{50}}
\put(1462,1932){\ellipse{50}{50}}
\put(1502,1837){\blacken\ellipse{50}{50}}
\put(1502,1837){\ellipse{50}{50}}
\put(1617,1777){\blacken\ellipse{50}{50}}
\put(1617,1777){\ellipse{50}{50}}
\put(1712,1777){\blacken\ellipse{50}{50}}
\put(1712,1777){\ellipse{50}{50}}
\put(1812,1897){\blacken\ellipse{50}{50}}
\put(1812,1897){\ellipse{50}{50}}
\put(1837,1987){\blacken\ellipse{50}{50}}
\put(1837,1987){\ellipse{50}{50}}
\put(2622,1872){\blacken\ellipse{50}{50}}
\put(2622,1872){\ellipse{50}{50}}
\put(2597,1777){\blacken\ellipse{50}{50}}
\put(2597,1777){\ellipse{50}{50}}
\put(2692,1657){\blacken\ellipse{50}{50}}
\put(2692,1657){\ellipse{50}{50}}
\put(2782,1667){\blacken\ellipse{50}{50}}
\put(2782,1667){\ellipse{50}{50}}
\put(2842,1867){\blacken\ellipse{50}{50}}
\put(2842,1867){\ellipse{50}{50}}
\put(2782,1947){\blacken\ellipse{50}{50}}
\put(2782,1947){\ellipse{50}{50}}
\put(2777,1257){\blacken\ellipse{50}{50}}
\put(2777,1257){\ellipse{50}{50}}
\put(2867,1277){\blacken\ellipse{50}{50}}
\put(2867,1277){\ellipse{50}{50}}
\put(2577,1107){\blacken\ellipse{50}{50}}
\put(2577,1107){\ellipse{50}{50}}
\put(2577,1017){\blacken\ellipse{50}{50}}
\put(2577,1017){\ellipse{50}{50}}
\put(2667,972){\blacken\ellipse{50}{50}}
\put(2667,972){\ellipse{50}{50}}
\put(2732,887){\blacken\ellipse{50}{50}}
\put(2732,887){\ellipse{50}{50}}
\put(2937,947){\blacken\ellipse{50}{50}}
\put(2937,947){\ellipse{50}{50}}
\put(3007,1002){\blacken\ellipse{50}{50}}
\put(3007,1002){\ellipse{50}{50}}
\put(3042,1122){\blacken\ellipse{50}{50}}
\put(3042,1122){\ellipse{50}{50}}
\put(2992,1197){\blacken\ellipse{50}{50}}
\put(2992,1197){\ellipse{50}{50}}
\put(2442,547){\blacken\ellipse{50}{50}}
\put(2442,547){\ellipse{50}{50}}
\put(2392,472){\blacken\ellipse{50}{50}}
\put(2392,472){\ellipse{50}{50}}
\put(2267,422){\blacken\ellipse{50}{50}}
\put(2267,422){\ellipse{50}{50}}
\put(2172,452){\blacken\ellipse{50}{50}}
\put(2172,452){\ellipse{50}{50}}
\put(2117,537){\blacken\ellipse{50}{50}}
\put(2117,537){\ellipse{50}{50}}
\put(2092,632){\blacken\ellipse{50}{50}}
\put(2092,632){\ellipse{50}{50}}
\put(2422,682){\blacken\ellipse{50}{50}}
\put(2422,682){\ellipse{50}{50}}
\put(2367,747){\blacken\ellipse{50}{50}}
\put(2367,747){\ellipse{50}{50}}
\put(1577,532){\blacken\ellipse{50}{50}}
\put(1577,532){\ellipse{50}{50}}
\put(1607,427){\blacken\ellipse{50}{50}}
\put(1607,427){\ellipse{50}{50}}
\put(1457,332){\blacken\ellipse{50}{50}}
\put(1457,332){\ellipse{50}{50}}
\put(1372,397){\blacken\ellipse{50}{50}}
\put(1372,397){\ellipse{50}{50}}
\put(1327,487){\blacken\ellipse{50}{50}}
\put(1327,487){\ellipse{50}{50}}
\put(1372,582){\blacken\ellipse{50}{50}}
\put(1372,582){\ellipse{50}{50}}
\put(1797,1117){\blacken\ellipse{50}{50}}
\put(1797,1117){\ellipse{50}{50}}
\put(1797,1017){\blacken\ellipse{50}{50}}
\put(1797,1017){\ellipse{50}{50}}
\put(1617,1227){\blacken\ellipse{50}{50}}
\put(1617,1227){\ellipse{50}{50}}
\put(1707,1227){\blacken\ellipse{50}{50}}
\put(1707,1227){\ellipse{50}{50}}
\put(1462,1202){\blacken\ellipse{50}{50}}
\put(1462,1202){\ellipse{50}{50}}
\put(1412,1117){\blacken\ellipse{50}{50}}
\put(1412,1117){\ellipse{50}{50}}
\put(1447,1037){\blacken\ellipse{50}{50}}
\put(1447,1037){\ellipse{50}{50}}
\put(1482,947){\blacken\ellipse{50}{50}}
\put(1482,947){\ellipse{50}{50}}
\put(977,727){\blacken\ellipse{50}{50}}
\put(977,727){\ellipse{50}{50}}
\put(932,632){\blacken\ellipse{50}{50}}
\put(932,632){\ellipse{50}{50}}
\put(962,887){\blacken\ellipse{50}{50}}
\put(962,887){\ellipse{50}{50}}
\put(987,802){\blacken\ellipse{50}{50}}
\put(987,802){\ellipse{50}{50}}
\put(692,872){\blacken\ellipse{50}{50}}
\put(692,872){\ellipse{50}{50}}
\put(782,897){\blacken\ellipse{50}{50}}
\put(782,897){\ellipse{50}{50}}
\put(612,757){\blacken\ellipse{50}{50}}
\put(612,757){\ellipse{50}{50}}
\put(657,667){\blacken\ellipse{50}{50}}
\put(657,667){\ellipse{50}{50}}
\put(587,1342){\blacken\ellipse{50}{50}}
\put(587,1342){\ellipse{50}{50}}
\put(692,1367){\blacken\ellipse{50}{50}}
\put(692,1367){\ellipse{50}{50}}
\put(772,1397){\blacken\ellipse{50}{50}}
\put(772,1397){\ellipse{50}{50}}
\put(807,1477){\blacken\ellipse{50}{50}}
\put(807,1477){\ellipse{50}{50}}
\put(807,1547){\blacken\ellipse{50}{50}}
\put(807,1547){\ellipse{50}{50}}
\put(767,1632){\blacken\ellipse{50}{50}}
\put(767,1632){\ellipse{50}{50}}
\put(712,1667){\blacken\ellipse{50}{50}}
\put(712,1667){\ellipse{50}{50}}
\put(622,1692){\blacken\ellipse{50}{50}}
\put(622,1692){\ellipse{50}{50}}
\put(522,1627){\blacken\ellipse{50}{50}}
\put(522,1627){\ellipse{50}{50}}
\put(472,1557){\blacken\ellipse{50}{50}}
\put(472,1557){\ellipse{50}{50}}
\path(12,462)(762,762)(1512,462)
	(2262,612)(2862,1062)(3612,462)
\path(12,1962)(612,1512)(762,762)
	(1662,1062)(2862,1062)(2712,1812)
	(1662,1962)(1662,1062)(612,1512)
	(1662,1962)(1512,2412)
\path(612,1512)(912,2412)
\path(2712,1812)(3462,2412)
\path(2862,1062)(3612,1512)
\path(2262,612)(3162,12)
\path(2262,612)(2112,12)
\path(1512,462)(1062,12)
\end{picture}
}